\documentclass{amsart}

\usepackage{amsmath,amssymb,amsthm,tikz,mathtools}
\usepackage{mathrsfs}
\usepackage{booktabs}
\usepackage{float}
\usepackage{needspace}
\usepackage[hidelinks]{hyperref}
\usepackage[numbers]{natbib}

\newtheorem{theorem}{Theorem}[section]
\newtheorem{lemma}[theorem]{Lemma}
\newtheorem{proposition}[theorem]{Proposition}
\newtheorem{corollary}[theorem]{Corollary}
\numberwithin{equation}{section}
\theoremstyle{remark}
\newtheorem{note}{Note}[section]

\theoremstyle{definition}
\newtheorem{definition}[theorem]{Definition}
\newtheorem{example}[theorem]{Example}
\newcommand{\norm}[1]{\left\lVert #1\right\rVert}
\newcommand{\transp}{\mathsf{T}}

\title[Veronese Classification in Complex Hyperquadrics]
{Veronese Classification of Constant-Curvature Minimal Two-Spheres
in Complex Hyperquadrics}

\author{Zonghao Jie}
\address{School of Mathematical Sciences,
University of Chinese Academy of Sciences}
\email{jiezonghao23@mails.ucas.ac.cn}
\author{Xiaoxiang Jiao}
\address{School of Mathematical Sciences,
University of Chinese Academy of Sciences}
\email{xxjiao@ucas.ac.cn}
\hypersetup{
  pdftitle={Veronese Classification of Constant-Curvature Minimal Two-Spheres in Complex Hyperquadrics},
  pdfauthor={Zonghao Jie and Xiaoxiang Jiao}
}



\begin{document}

\begin{abstract}
We prove that every linearly full, irreducible conformal minimal
immersion of constant Gaussian curvature and finite isotropy order
from the two-sphere into a complex hyperquadric comes from a member of
a Veronese sequence. The result holds at every finite isotropy order,
without an unramifiedness or additional Frenet-norm polynomiality
assumption. Both the associated projective harmonic sequence and the
target harmonic sequence are totally unramified as consequences. If the
replacement depth is \(\ell\) and the holomorphic seed spans
\(\mathbb{CP}^m\), the curvature is
\(K=2/((2\ell+1)m-2\ell^2)\) in the Grassmannian normalization.
Real Gram data characterize the admissible embeddings and give an
algebraic converse. The proof uses holomorphic Pfaffian norm identities,
degree estimates that retain all ramification divisors, and a valuation
argument that identifies the retained lines with seed Frenet lines.
\end{abstract}

\maketitle
\vspace*{-3pt}

\section{Introduction}

Let \(\varphi:\mathbb S^2\to G(2,n;\mathbb R)\) be a linearly full,
irreducible conformal minimal immersion of finite isotropy order
\(r=2s+1\). We fix an isotropic complex-line lift through the isometric
double cover \(Q_{n-2}\to G(2,n;\mathbb R)\). The replacement construction
of Bahy-El-Dien and Wood associates to \(\varphi\) a holomorphic seed
and a finite tower of alternating forward and backward replacements.
We prove that constant Gaussian curvature identifies every retained
line with a Frenet line of a Veronese seed. No unramifiedness assumption
is required. The proof combines local replacement identities with
global degree estimates that retain all ramification divisors.

\noindent\textbf{Relation with earlier work.}
Homogeneous minimal two-spheres in complex hyperquadrics were classified
by Peng, Wang, and Xu \citep{peng2016homogeneous}. For the constant-curvature
problem without a homogeneity assumption, Jiao and M.~Li
\citep[Proposition~4.3 and Theorem~4.6(iii)]{jiao2017classification}
treat the linearly full, irreducible, totally unramified case with finite
isotropy order \(r\geq n-5\). Their depth-one classification has
\(m=n-1\) and curvature \(2/(3n-5)\). H.~Li and Jiao
\citep[Propositions~3--4 and Theorem~1]{li2023conformal} consider
\(r=n-6\). At depth two, with total unramifiedness and the additional
polynomiality of \(\norm{f_2^{(n-1)}}^2\), they obtain the second Veronese
real layer and curvature \(2/(5n-13)\). The Pfaffian gradings below
extend their low-order holomorphic quantities to arbitrary depth,
without an additional polynomiality assumption on a Frenet norm.

Low-dimensional classifications in \(Q_2\), \(Q_3\), \(Q_4\), and
\(Q_5\) are given in
\citep{jiao2013q2,li2016q3,jiao2021q4,jiao2019q5}, with different
ramification and isotropy hypotheses. We compare their classified
families in \S\ref{subsec:known-examples}. In particular, the
nonhomogeneous \(Q_4\) example of Jiao, M.~Li, and H.~Li is reducible;
their irreducible totally unramified examples are strongly isotropic.

Chi, Xie, and Xu \citep{chi2021structure} study constant-curvature spheres
that are minimal both in the quadric and in its ambient projective space,
and describe their real extrinsic data by symmetric matrices. Their
projective-minimality hypothesis is not automatic from quadric
minimality alone. Here it follows from the vanishing of every replacement
deviation. The resulting real Gram description is consistent with that
work and is not a claim that all these immersions are homogeneous or
that their full moduli spaces have been explicitly solved.

Roundness and Veronese rigidity for primitive flag lifts of holomorphic
Grassmannian sequences are also obtained by Pacheco and Rehman
\citep{pacheco2026primitive}. Their positive invariant metric weights
all successive connections. In the present setting, before faithfulness
is proved there is no single projective Frenet flag of retained lines,
and the metric controls only two adjacent Pfaffian prefixes. The mixed
kernel and its source equation provide the additional backward
propagation. Our external replacement and rationality inputs are the
construction theorems of \citep{burstall1986construction,bahy1989explicit};
they do not presuppose unramifiedness or a constant-curvature
classification.

\noindent\textbf{The mechanism.}
The argument has three main steps.
\begin{enumerate}
\item Graph connection and contact identities give holomorphic
Pfaffian gradings \(\Psi_{j,\nu}\) and the raw norm identity
\[
 P_j:=\prod_{i=0}^j\norm{v_i}^2
 =\sum_\nu\norm{\Psi_{j,\nu}}^2.
\]
These are initially rational data, not globally polynomial lifts in
a common regular gauge. Independent saturation gives polynomial
vectors of actual degrees \(D_j\), whose curvature forms equal the
smooth forward energies of the genuine retained lines.
\item Finite isotropy gives a nonzero holomorphic bilinear square at
the middle seed line. Together with the highest-row projection of a
retained line onto its seed Frenet line, its effective divisor yields
\(D_j-D_{j-1}\geq2\). All ramification and projection zeros are retained
in this calculation; the standard Veronese degrees are not yet known.
\item Constant curvature rounds the final two saturated norms. Exact
norm and source equations in the original rational gauge force any
nonzero mixed kernel to descend to an impossible valuation at the
bottom of the tower. Hence every deviation vanishes, the seed norm
rounds, and coefficient comparison identifies the Veronese curve.
Only afterwards do we deduce the standard degrees and independently
verify total unramifiedness of the target.
\end{enumerate}

Here \(\ell\) is the actual depth of the noncancelling replacement
construction in Lemma~\ref{lem:bew-replacement-series}; an immediately
inverse pair is not counted. This is the canonical construction depth,
not an additional generic-position condition. The integer \(m\) is
the dimension of the projective span of the seed and need not equal
\(n-1\). Curvature is measured in
\eqref{eq:quadric-metric-normalization}.

\begin{theorem}[Veronese classification]\label{thm:main-veronese-classification}
Let \(\varphi:\mathbb S^2\to G(2,n;\mathbb R)\) be a linearly full,
irreducible conformal minimal immersion of constant Gaussian curvature
and finite isotropy order \(r=2s+1\), \(s\geq0\). Fix an isotropic
complex-line lift to \(Q_{n-2}\). Let \(\ell\geq1\) be its replacement
depth and \(\mathbb{CP}^m\) the projective span of its holomorphic seed.
Then the seed has degree \(m\), the tower is totally faithful, and
both the seed's projective harmonic sequence and the target's two-sided
harmonic sequence are totally unramified.

There is a spherical affine coordinate \(z\), a complex linear isometric
embedding \(U:\mathbb C^{m+1}\hookrightarrow\mathbb C^n\), and constants
\(C_j>0\) such that, in the monic normalization,
\begin{gather}
 v_j=f_j\quad(0\leq j\leq\ell),\qquad
 f_0=\sqrt{C_0}\,UV_0^{(m)},\label{eq:main-faithfulness}\\
 P_j=C_j(1+|z|^2)^{d_j},\qquad
 d_j=(j+1)(m-j),\quad0\leq j\leq\ell,\label{eq:main-prefix-degrees}\\
 \underline\varphi
 =\operatorname{span}(UV_\ell^{(m)})\oplus
  \operatorname{span}(\overline{UV_\ell^{(m)}}).
 \label{eq:main-real-layer}
\end{gather}
Here \(V_j^{(m)}\) is the standard Veronese Frenet vector.
In particular, the isotropic complex-line lift is minimal in
\(\mathbb{CP}^{n-1}\). The curvature is
\begin{equation}\label{eq:main-curvature}
 \boxed{K=\frac{2}{(2\ell+1)m-2\ell^2}.}
\end{equation}
Conversely, an isometric embedding satisfying the finite-isotropy
conditions \eqref{eq:realgram-isotropy} and the real-fullness condition
\eqref{eq:realgram-fullness} gives such an immersion, as specified in
Corollary~\ref{cor:veronese-converse}.
\end{theorem}

The admissible range supplied by finite isotropy is
\begin{equation}\label{eq:intro-depth-bound}
 2(\ell+s)+2\leq m\leq n-1.
\end{equation}
The theorem concerns immersions of finite isotropy order and irreducible
towers; it does not assert a classification of branched target maps,
reducible maps, or strongly isotropic spheres. No polynomiality or
ambient projective-minimality assumption is added.

\noindent\textbf{Organization.}
Section~\ref{sec:preliminaries} fixes the metric and the distinction
between saturated bundles and derivative lifts.
Section~\ref{sec:basic-setup} gives the replacement and short-window
graph calculus. Section~\ref{sec:pfaffian-data} constructs the rational
Pfaffian data and their independent saturations.
Section~\ref{sec:degree-growth} proves degree growth with effective
divisors, and Section~\ref{sec:rigidity} proves raw-gauge backward
rigidity. Section~\ref{sec:automatic-unramifiedness} obtains the
Veronese representation and both unramifiedness conclusions.
Section~\ref{sec:real-gram} gives the converse, examples, and the
conclusion of the main theorem. Appendix~\ref{app:graph-connection}
contains the full graph induction. Appendix~\ref{app:bordered-pfaffian}
supplies the full bordered-Pfaffian calculation used in the mixed kernel.

    \section{Preliminaries}\label{sec:preliminaries}
     \subsection{The Relationship between
     \texorpdfstring{\(Q_{n-2}\) and \(G(2,n;\mathbb{R})\)}
     {Q(n-2) and G(2,n;R)}}
    
    We first define the Hermitian inner product and the symmetric complex
    bilinear form by
    \begin{equation}\label{eq:inner-product}
    	H(z,w)=z^{\transp}\overline{w},
    	\qquad
    	b(z,w)=z^{\transp}w,
    \end{equation}
    where \(z,w\in\mathbb{C}^n\).
    
    Let \(\Bbbk\) be a field and let \(\Bbbk^n\) be the \(n\)-dimensional
    \(\Bbbk\)-vector space. Denote by \(G(k,n;\Bbbk)\) the set of
    \(k\)-dimensional linear subspaces of \(\Bbbk^n\), where \(0\leq k\leq n\).
    This set is called the \(\Bbbk\)-Grassmannian; in the applications below,
    \(\Bbbk\) is either \(\mathbb{R}\) or \(\mathbb{C}\).
    
    Let
    \begin{equation}\label{eq:quadric}
    	Q_{n-2}
    	=\left\{[z]=[z_1,\dots,z_n]\in\mathbb{CP}^{n-1}
    	\;\middle|\;
    	b(z,z)=\sum_{k=1}^n z_k^2=0\right\}
    \end{equation}
    denote the complex quadric hypersurface. Let \(g_{\mathrm{FS}}\)
    be the Fubini--Study metric of holomorphic sectional curvature
    \(4\), and let
    \(\sigma:Q_{n-2}\to\mathbb{CP}^{n-1}\) be the inclusion. The map
    \begin{equation}\label{eq:quadric-covering}
    	f:Q_{n-2}\to G(2,n;\mathbb{R}),
    	\qquad
    	[Z]\mapsto\operatorname{span}_{\mathbb{R}}
    	\{\operatorname{Re}Z,\operatorname{Im}Z\},
    \end{equation}
    is a well-defined double cover \citep{wolfson1986harmonic}.
    Throughout the paper, the real Grassmannian carries the metric
    induced by \eqref{eq:grassmannian-metric}. In this normalization,
    we equip the quadric with
    \begin{equation}\label{eq:quadric-metric-normalization}
    	g_Q=f^*g_{G(2,n;\mathbb{R})}
    	=2\sigma^*g_{\mathrm{FS}}.
    \end{equation}
    Indeed, an isotropic line with orthogonal projection \(\Pi\)
    corresponds to the real projection \(\Pi+\overline\Pi\);
    substituting this into \eqref{eq:grassmannian-metric} gives the
    factor \(2\). Thus \(f\) is a local isometry for \(g_Q\), whereas
    \(\sigma\) scales squared lengths by \(1/2\). The curvature \(K\)
    in the main theorem refers to \(g_Q\), or equivalently to the
    Grassmannian metric. With the induced Fubini--Study metric on the
    quadric, the same immersion has Gaussian curvature \(2K\).
    Since \(\mathbb{S}^2\) is simply connected, a map
    \(\mathbb{S}^2\to G(2,n;\mathbb{R})\) has a lift to \(Q_{n-2}\),
    unique once a lift at one point is chosen. The complexified real
    bundle of a map with isotropic line lift \(L\) is
    \(L\oplus\overline L\); this rank-two bundle is distinguished from
    the chosen rank-one lift throughout.

    We identify \(\mathbb{S}^2\) biholomorphically with
    \(\mathbb{CP}^1\). A spherical affine coordinate \(z\) means the
    standard coordinate on the complement of one point under such an
    identification; near infinity we use \(\zeta=1/z\).
    
    \subsection{Harmonic Sequences in
    \texorpdfstring{\(G(k,n;\mathbb{C})\)}{G(k,n;C)}}
    Let \(M\) be a Riemann surface and let
    \(\varphi:M\to G(k,n;\mathbb{C})\). The map \(\varphi\) corresponds to a
    rank-\(k\) subbundle of the trivial bundle
    \(\underline{\mathbb{C}^n}=M\times\mathbb{C}^n\), denoted by
    \(\underline{\varphi}\). Its fiber over \(x\) is
    \(\underline{\varphi}_x = \varphi(x)\). The map \(\varphi\) is
    holomorphic if and only if \(\underline\varphi\) is a holomorphic
    subbundle of the trivial holomorphic bundle, equivalently if its
    \(\bar\partial\)-second fundamental form vanishes. If \(\varphi\) is harmonic,
    \(\underline{\varphi}\) is called a harmonic subbundle of \(\underline{\mathbb{C}^n}\). If
    \(\underline{\varphi}\) is not contained in any proper trivial subbundle of
    \(\underline{\mathbb{C}^n}\) of the form
    \(\underline{\mathbb{C}^m}\) with \(m<n\), then \(\varphi\) is said to be
    linearly full.
    
    Let \(z\) be a local complex coordinate on \(M\), and write
    \begin{equation}\label{eq:dolbeault-operators}
    	\begin{aligned}
    	\partial &= \frac{\partial}{\partial z}, &
    	\bar\partial &= \frac{\partial}{\partial\overline{z}}.
    \end{aligned}
    \end{equation}
    Let \(\underline{\varphi}\) and \(\underline{\psi}\) be mutually orthogonal
    subbundles of the trivial bundle \(\underline{\mathbb{C}^n}\). We define
    \begin{equation}\label{eq:second-fundamental-forms}
    	\begin{aligned}
    	A'_{\underline{\varphi},\underline{\psi}}:
    	\underline{\varphi}\to\underline{\psi},\quad
    	&v\mapsto\pi_{\underline{\psi}}(\partial v),\\
    	A''_{\underline{\varphi},\underline{\psi}}:
    	\underline{\varphi}\to\underline{\psi},\quad
    	&v\mapsto\pi_{\underline{\psi}}(\bar\partial v),
    \end{aligned}
    \end{equation}
    These maps are the \(\partial'\)- and \(\partial''\)-second fundamental
    forms, respectively, where \(v\) is a local smooth section of
    \(\underline{\varphi}\). We write
    \(A'_{\underline{\varphi}}=A'_{\underline{\varphi},\underline{\varphi}^\perp},A''_{\underline{\varphi}}=A''_{\underline{\varphi},\underline{\varphi}^\perp}\).
    
    \begin{lemma}[Holomorphicity of the harmonic differentials]
    	\label{lem:harmonic-connection-holomorphic}
    	\citep{burstall1986construction, wolfson1988harmonic}
    	Suppose that \(\varphi:M\to G(k,n;\mathbb{C})\) is a harmonic map.
    	With the induced bundle structures,
    	\(A'_{\underline{\varphi}}\) is holomorphic and
    	\(A''_{\underline{\varphi}}\) is antiholomorphic. The maps
    	\(\varphi_1\) and
    	\(\varphi_{-1}\) are also harmonic maps, where
    	\(\underline{\varphi_1} = \operatorname{Im}A'_{\underline{\varphi}},\underline{\varphi_{-1}} = \operatorname{Im}A''_{\underline{\varphi}}\).
    	Here and below, images are extended as smooth saturated
    	subbundles across isolated rank drops of the differentials.
    \end{lemma}
    
    Thus every harmonic map \(\varphi:M\to G(k,n;\mathbb{C})\) gives rise to
    the following two harmonic sequences:
    
    \begin{equation}\label{eq:harmonic-sequence-forward}
    	\underline{\varphi}
    	=\underline{\varphi_0}
    	\stackrel{\partial'}{\rightarrow}\underline{\varphi_1}
    	\stackrel{\partial'}{\rightarrow}\cdots
    	\stackrel{\partial'}{\rightarrow}\underline{\varphi_\alpha}
    	\stackrel{\partial'}{\rightarrow}\cdots,
    \end{equation}
    \begin{equation}\label{eq:harmonic-sequence-backward}
    	\underline{\varphi}
    	=\underline{\varphi_0}
    	\stackrel{\partial''}{\rightarrow}\underline{\varphi_{-1}}
    	\stackrel{\partial''}{\rightarrow}\cdots
    	\stackrel{\partial''}{\rightarrow}\underline{\varphi_{-\alpha}}
    	\stackrel{\partial''}{\rightarrow}\cdots.
    \end{equation}
    where
    \(\underline{\varphi_\alpha}=\partial'\underline{\varphi_{\alpha -1}}=\operatorname{Im}A'_{\underline{\varphi_{\alpha -1}}}\),
    \(\underline{\varphi_{-\alpha}} = \partial''\underline{\varphi_{-\alpha +1}} = \operatorname{Im}A''_{\underline{\varphi_{{-\alpha} +1}}}\),
    for \(\alpha>0\). The subbundles in these sequences are called Gauss
    bundles.

    \begin{definition}[Gauss bundles]
    	\label{def:gauss-bundles}
    	\citep{burstall1986construction,wolfson1988harmonic}
    	Let \(\varphi:M\to G(k,n;\mathbb{C})\) be harmonic, and identify
    	\(\varphi\) with its associated rank-\(k\) subbundle
    	\(\underline{\varphi}\). For every integer \(j\), define the
    	\(j\)-th Gauss bundle by
    	\begin{equation}\label{eq:gauss-bundles}
    		G^{(j)}(\varphi):=\underline{\varphi_j},
    	\end{equation}
    	where \(\underline{\varphi_0}=\underline{\varphi}\) and, for
    	\(\alpha>0\),
    	\[
    		\underline{\varphi_\alpha}
    		=\partial'\underline{\varphi_{\alpha-1}},
    		\qquad
    		\underline{\varphi_{-\alpha}}
    		=\partial''\underline{\varphi_{-\alpha+1}}.
    	\]
    	Thus \(G^{(0)}(\varphi)=\underline{\varphi}\), while
    	\(G^{(j)}(\varphi)\) is the \(j\)-th forward Gauss bundle when
    	\(j>0\), and the \((-j)\)-th backward Gauss bundle when \(j<0\).
    	We also write \(G'(\varphi)=G^{(1)}(\varphi)\) and
    	\(G''(\varphi)=G^{(-1)}(\varphi)\).
    \end{definition}
    
    Lemma~\ref{lem:harmonic-sequence-properties} collects the
    properties of these sequences.
    
    \begin{lemma}[Properties of harmonic sequences]
    	\label{lem:harmonic-sequence-properties}
    	\citep{burstall1986construction}
    	Let \(\varphi\) be harmonic.
    	
    	(1)
    	\(\partial''(\partial' \underline{\varphi_\alpha})\subset \underline{\varphi_\alpha},\partial'(\partial'' \underline{\varphi_\alpha})\subset \underline{\varphi_\alpha}\).
    	In particular, if \(\operatorname{rank}\varphi_\alpha=\operatorname{rank}\varphi_{\alpha+1}\), then \(\partial''(\partial' \underline{\varphi_\alpha})=\underline{\varphi_\alpha}\)
    	and \(\partial'(\partial'' \underline{\varphi_{\alpha +1}})=\underline{\varphi_{\alpha +1}}\).
    	
    	(2) If \(M=\mathbb{S}^2\), then the harmonic sequences terminate after finitely many steps; that is,
    	there exist \(\alpha_1,\alpha_2>0\) such that \(\varphi_{\alpha_1}=\varphi_{-\alpha_2}=0\).
    \end{lemma}
    
    \begin{definition}\label{def:harmonic-sequence-reducibility-isotropy}
    	\citep{burstall1986construction}
    	For a harmonic map \(\varphi\), define the following:
    	
    	(1) If \(\operatorname{rank}\underline{\varphi_i}>\operatorname{rank}\underline{\varphi_{i+1}}\) (respectively,
    	\(\operatorname{rank}\underline{\varphi_i}>\operatorname{rank}\underline{\varphi_{i-1}}\)), then \(\varphi_i\) is
    	called \(\partial'-\)reducible (respectively, \(\partial''-\)reducible); otherwise, it is called
    	\(\partial'-\)irreducible (respectively, \(\partial''-\)irreducible).
    	
    	(2) Set
    	\[
    	 r=\sup\{q\in\mathbb{Z}_{\geq0}:
    	 \underline{\varphi_i}\perp\underline{\varphi}
    	 \text{ for every }1\leq i\leq q\}.
    	\]
    	This is the isotropy order of \(\varphi\). For finite \(r\),
    	\(\underline{\varphi_{r+1}}\) is not orthogonal to
    	\(\underline{\varphi}\). If \(r=\infty\), then \(\varphi\) is
    	called strongly isotropic.

    	Equivalently, \(r\) is the largest integer such that
    	\(G^{(i)}(\varphi)\perp G^{(j)}(\varphi)\) whenever
    	\(0<|i-j|\leq r\); see \citep[\S1(C)]{bahy1989explicit}.
    	
    	(3) If \(\varphi\) has finite isotropy order, define
    	\begin{equation}
    		c'_r(\varphi)=A'_{\underline{\varphi_r},\underline{\varphi}} \circ A'_{\underline{\varphi_{r-1}},\underline{\varphi_{r}}} \circ \cdots \circ A'_{\underline{\varphi},\underline{\varphi_{1}}},
    	\end{equation}
    	This composition is called the first \(\partial'\)-closing map. The first
    	\(\partial''\)-closing map is defined analogously.
    \end{definition}
    
    Lemma~\ref{lem:real-grassmannian-properties} records the properties
    of harmonic maps into the real Grassmannian
    \(G(2,n;\mathbb{R})\).
    
    \begin{lemma}[Properties of harmonic maps into real Grassmannians]
    	\label{lem:real-grassmannian-properties}
    	\citep{bahy1989explicit}
    	Let \(\varphi:\mathbb{S}^2\to G(2,n;\mathbb{R})\) be harmonic.
    	
    	(1) \(\varphi\) is \(\partial'-\)reducible if and only if it is \(\partial''-\)reducible.
    	
    	(2) If the isotropy order \(r\) of \(\varphi\) is finite, then \(r\) is odd. More precisely, if
    	\(r=2s+1\), then
    	\(0\leq s\leq\left\lfloor\frac{n-3}{2}\right\rfloor\).
    	
    	(3) If \(r<\infty\), then
    	\(\underline{\varphi} = \underline{V} \oplus \underline{\overline{V}}\), where
    	\(\overline{\underline{V}}\) is the image of the first
    	\(\partial'\)-closing map, namely
    	\(\operatorname{Im}c'_r(\varphi)\).
    \end{lemma}
    
    A conformal immersion from a Riemann surface is minimal precisely
    when it is harmonic \citep{eells1964harmonic}.
    
    \subsection{Geometry in
    \texorpdfstring{\(G(k,n;\mathbb{C})\)}{G(k,n;C)}}
    \label{sec:geometry-complex-grassmannian}
    
    Given
    \(\varphi:\mathbb{S}^2\to G(k,n;\mathbb{C})\), choose a local orthonormal
    frame \(\{e_1,\dots,e_k\}\) for \(\underline{\varphi}\), and set
    \(W=(e_1,\dots,e_k)\). The map is then represented by
    \(\varphi=WW^*\), a matrix that is independent of the chosen orthonormal
    frame. This identifies \(G(k,n;\mathbb{C})\) with the set of rank-\(k\)
    orthogonal projections on \(\mathbb{C}^n\):
    \begin{equation}\label{eq:complex-grassmannian-projections}
    	G(k,n;\mathbb{C})
    	=\left\{\varphi\in M_n(\mathbb{C})
    	\;\middle|\;
    	\varphi^2=\varphi,
    	\;\varphi^*=\varphi,
    	\;\operatorname{rank}\varphi=k\right\}.
    \end{equation}
    The map \(\iota(\varphi)=2\varphi-I\) embeds
    \(G(k,n;\mathbb{C})\) totally geodesically into \(U(n)\).
    
    Equip \(U(n)\) with its Maurer--Cartan form \(\omega\) and the metric
    \(\mathrm{ds}^2_{U(n)}=\frac{1}{8}\operatorname{tr}\omega\omega^*\).
    The metric pulled back by \(\iota\) is
    \begin{equation}\label{eq:grassmannian-metric}
    	\mathrm{ds}^2_{G(k,n;\mathbb{C})}=\frac{1}{2}\operatorname{tr}\mathrm{d}\varphi\mathrm{d}\varphi.
    \end{equation}
    For a conformal map, comparison of coefficients gives the induced
    metric
    \begin{equation}\label{eq:conformal-grassmannian-metric}
    	\mathrm{ds}^2=\operatorname{tr}
    	\bigl((\partial\varphi)(\bar\partial\varphi)\bigr)
    	\,\mathrm{d}z\,\mathrm{d}\overline z.
    \end{equation}
    Let \(\varphi:\mathbb{S}^2\to G(k,n;\mathbb{C})\) be a linearly full
    harmonic map whose harmonic sequence is
    
    \begin{equation}
    	\underline{\varphi_0}\stackrel{\partial'}{\rightarrow}\cdots
    	\stackrel{\partial'}{\rightarrow}\underline{\varphi_{j-1}}
    	\stackrel{\partial'}{\rightarrow}\underline{\varphi}=\underline{\varphi_j}
    	\stackrel{\partial'}{\rightarrow}\underline{\varphi_{j+1}}
    	\stackrel{\partial'}{\rightarrow}\cdots
    	\stackrel{\partial'}{\rightarrow}\underline{\varphi_{j_0}}
    	\stackrel{\partial'}{\rightarrow}0.
    \end{equation}
    Write \(\operatorname{rank}\varphi_i=k_i\), and choose a local
    orthonormal frame
    \begin{equation}
    	\{ e_1^{(i)},e_2^{(i)},\cdots,e_{k_i}^{(i)} \}
    \end{equation}
    for \(\underline{\varphi_i}\). Set
    \begin{equation}
    	W_i=(e_1^{(i)},e_2^{(i)},\cdots,e_{k_i}^{(i)}),\qquad i=0,1,\cdots , j_0.
    \end{equation}
    Then
    \begin{equation}
    	\begin{aligned}
    	W_iW_i^*&=\varphi_i, & W_i^*W_i&=I_{k_i\times k_i},\\
    	W_i^*W_{i+1}&=0_{k_i\times k_{i+1}}, &
    	W_i^*W_{i-1}&=0_{k_i\times k_{i-1}}.
    \end{aligned}
    \end{equation}
    Write
    \(\partial W_i=W_{i+1}\Omega_i+W_i\Psi_i\), where \(\Omega_i\) is a
    \(k_{i+1}\times k_i\) matrix and \(\Psi_i\) is a \(k_i\times k_i\)
    matrix. A direct calculation gives
    \(\bar\partial W_i=-W_{i-1}\Omega_{i-1}^*-W_i\Psi_i^*\). The
    integrability conditions are
    
    \begin{equation}\label{eq:connection-integrrability}
    	\begin{cases}
    	\bar\partial\Omega_i = \Psi_{i+1}^*\Omega_i - \Omega_i\Psi_i^* \\
    	\bar\partial\Psi_i + \partial\Psi_i^* = \Omega_i^*\Omega_i + \Psi_i^*\Psi_i - \Omega_{i-1}\Omega_{i-1}^* - \Psi_i\Psi_i^*
    	\end{cases}.
    \end{equation}
    Set \(L_i=\operatorname{tr}(\Omega_i\Omega_i^*)\). The metric
    induced by \(\varphi_i\) is
    \begin{equation}\label{eq:connection-density-metric}
    	\mathrm{ds}_i^2=(L_{i-1}+L_i)\mathrm{d}z\mathrm{d}\overline{z}.
    \end{equation}
    The Laplace--Beltrami operator and Gaussian curvature associated with \(\varphi_i\) are
    \begin{equation}\label{eq:connection-gaussian-curvature}
    	\Delta_i=\frac{4}{L_{i-1}+L_i}\partial\bar\partial,
    	\qquad
    	K_i=-\frac{2}{L_{i-1}+L_i}
    	\partial\bar\partial\log(L_{i-1}+L_i).
    \end{equation}
    The nondegeneracy of second fundamental forms along a harmonic sequence
    plays an important role in the study of pseudoholomorphic curves; see,
    for example, \citep{jiao2003pseudoholomorphic}. Following the later work
    on real Grassmannians and complex hyperquadrics
    \citep{jiao2017classification,li2023conformal}, we use the following
    convention for unramified maps.
    
    \begin{definition}\label{def:unramified}
    	For a harmonic map \(\varphi=\varphi_0\), consider a nonterminal
    	forward connection
    	\(A'_{\underline{\varphi_\alpha}}:
    	\underline{\varphi_\alpha}\to
    	\underline{\varphi_{\alpha+1}}\), with \(\alpha\geq0\) and
    	\(k_{\alpha+1}=\operatorname{rank}\varphi_{\alpha+1}>0\).
    	It is \emph{unramified} if the tensor
    	density
    	\begin{equation}\label{eq:unramified-density}
    		\det\bigl(\Omega_\alpha\Omega_\alpha^*\bigr)
    	\bigl(\mathrm{d}z\,\mathrm{d}\overline{z}\bigr)^{k_{\alpha+1}}
    	\end{equation}
    	is nowhere zero on \(\mathbb{S}^2\). A nonterminal backward
    	connection from \(\varphi_{-\alpha}\) to
    	\(\varphi_{-\alpha-1}\), with \(\alpha\geq0\), is unramified
    	when the analogous density for its \(\bar\partial\)-connection
    	matrix is nowhere zero. We call \(\varphi\)
    	\emph{totally unramified} if every nonterminal connection in
    	both Gauss sequences is unramified, including the two
    	connections starting at \(\varphi_0\).
    \end{definition}
    
    Here \(\Omega_\alpha\) is a
    \(k_{\alpha+1}\times k_\alpha\) matrix, so
    \(\Omega_\alpha\Omega_\alpha^*\) is a square matrix of size
    \(k_{\alpha+1}\). Its determinant is unchanged by unitary changes of
    local frames, while the coordinate factor compensates for a change of
    complex coordinate. In particular, the nowhere-vanishing condition is
    equivalent to \(\Omega_\alpha\) having full row rank at every point.
    The product \(\mathrm{d}z\,\mathrm{d}\overline{z}\) is interpreted as a
    tensor density rather than a wedge product. At a rank-two to
    rank-one interface this requires a nonzero row, and at an
    equal-rank interface it requires a pointwise isomorphism. The
    condition concerns the actual differential, not just its saturated
    image. For a real map the backward conditions follow from the
    forward ones by conjugation.
    
    \subsection{Projective Frenet and Veronese Sequences}
    
    Let \(\varphi:\mathbb{S}^2\to\mathbb{CP}^n\) be a linearly full
    harmonic map. Its
    harmonic sequence is
    
    \begin{equation}\label{eq:veronese-sequence}
    	0\stackrel{\partial''}{\leftarrow}\underline{\varphi_0}
    	\stackrel{\partial''}{\leftarrow}\cdots
    	\stackrel{\partial''}{\leftarrow}\underline{\varphi_{j-1}}
    	\stackrel{\partial''}{\leftarrow}\underline{\varphi}=\underline{\varphi_j}
    	\stackrel{\partial'}{\rightarrow}\underline{\varphi_{j+1}}
    	\stackrel{\partial'}{\rightarrow}\cdots
    	\stackrel{\partial'}{\rightarrow}\underline{\varphi_{j_0}}
    	\stackrel{\partial'}{\rightarrow}0,
    \end{equation}
    The properties of the sequence
    \eqref{eq:veronese-sequence} are summarized in
    Lemma~\ref{lem:veronese-frenet}.
    
    \begin{lemma}[Projective Frenet frame]
    	\label{lem:veronese-frenet}
    	\citep{bolton1988conformal}
    	Let \(f_0^{(n)}\) be a holomorphic section of
    	\(\underline{\varphi_0}\), and define
    	\(f_{q+1}^{(n)}
    	=\pi_{\underline{\varphi_q}^{\perp}}(\partial f_q^{(n)})\), where
    	\(q=0,1,\cdots,n-1\).
    	
    	(1) The family \(\{f_q^{(n)}\}_{q=0}^n\) is obtained from
    	\(\{\partial^q f_0^{(n)}\}_{q=0}^n\) by Gram--Schmidt
    	orthogonalization.
    	
    	(2) The vectors \(\{f_q^{(n)}\}_{q=0}^n\) are smooth sections of
    	\(\{\varphi_q\}_{q=0}^n\). Moreover, \(j_0=n\).
    	
    	(3) Set
    	\begin{equation}\label{eq:lambda-definition}
    	 \lambda_i^{(n)}
    	 =\frac{\norm{f_{i+1}^{(n)}}^2}{\norm{f_i^{(n)}}^2},
    	 \qquad 0\leq i<n,
    	\end{equation}
    	and \(\lambda_{-1}^{(n)}=\lambda_n^{(n)}=0\). These vectors satisfy
    	\begin{equation}\label{eq:veronese-frenet-derivative}
    		\partial f_{q}^{(n)}
    		=f_{q+1}^{(n)}
    		+\frac{\langle\partial f_{q}^{(n)},f_{q}^{(n)}\rangle}
    		{\lVert f_{q}^{(n)}\rVert^2}f_{q}^{(n)},
    		\qquad q=0,1,\dots,n-1,
    	\end{equation}
    	\begin{equation}\label{eq:veronese-bar-frenet-derivative}
    		\bar\partial f_{q}^{(n)}
    		=-\lambda_{q-1}^{(n)}f_{q-1}^{(n)},
    		\qquad q=1,2,\cdots,n.
    	\end{equation}
    \end{lemma}

    The quantities \(\lambda_i^{(n)}\) are the connection densities of
    the projective Frenet sequence. Their expression in terms of
    \(\beta_i\) and normalized frames appears after the local frame
    equations below.
    
    The identities in Lemma~\ref{lem:veronese-frenet} are first read
    on the common regular set. To include ramification points, write
    \(\mathcal L_i\) for the smooth saturated Frenet lines. Their
    induced holomorphic structures make the Frenet homomorphisms
    holomorphic \citep{wolfson1988harmonic}. In nowhere-zero
    induced-holomorphic frames \(e_i\), write
    \(\partial'(\partial_z\otimes e_i)=\beta_i e_{i+1}\),
    where \(\beta_i\) is holomorphic. The corresponding equations are
    \begin{align}
    \partial e_i
    &=\beta_i e_{i+1}
      +(\partial\log\norm{e_i}^2)e_i,
      \label{eq:frenet-forward}\\
    \bar\partial e_i
    &=-\overline{\beta_{i-1}}
      \frac{\norm{e_i}^2}{\norm{e_{i-1}}^2}e_{i-1},
      \label{eq:frenet-backward}\\
    \partial\overline{e_i}
    &=-\beta_{i-1}
      \frac{\norm{e_i}^2}{\norm{e_{i-1}}^2}\overline{e_{i-1}},
      \label{eq:frenet-conjugate-forward}\\
    \bar\partial\overline{e_i}
    &=\overline{\beta_i}\,\overline{e_{i+1}}
      +(\bar\partial\log\norm{e_i}^2)\overline{e_i}.
      \label{eq:frenet-conjugate-backward}
    \end{align}
    At an endpoint, the absent neighbouring term is omitted.
    In these frames,
    \begin{equation}\label{eq:lambda-in-frame}
     \lambda_i^{(n)}
     =|\beta_i|^2\frac{\norm{e_{i+1}}^2}{\norm{e_i}^2}.
    \end{equation}
    Thus a zero of a derivative lift is recorded by \(\beta_i\),
    rather than by a vanishing local frame. The \(i\)-th Frenet
    homomorphism is unramified when \(\beta_i\) has no zeros; the
    seed is totally unramified when this holds at every
    nonterminal index.

    Substituting
    \eqref{eq:veronese-frenet-derivative} and
    \eqref{eq:veronese-bar-frenet-derivative} into the formulas in
    \S\ref{sec:geometry-complex-grassmannian} gives
    
    \begin{equation}\label{eq:veronese-density}
    	L_q=\lambda_q^{(n)},
    	\qquad q=0,1,\dots,n-1,
    \end{equation}
    \begin{equation}\label{eq:veronese-toda}
    	\partial\bar\partial\log\lambda_q^{(n)}
    	=\lambda_{q+1}^{(n)}-2\lambda_q^{(n)}
    	+\lambda_{q-1}^{(n)},
    	\qquad q=0,1,\cdots,n-1,
    \end{equation}
    where the last term is omitted for \(q=0\) and the first
    term on the right is omitted for \(q=n-1\).
    
    Let
    \begin{equation}\label{eq:osculating-curve}
    	\begin{aligned}
    	\sigma_q^{(n)}=F_q^{(n)}
    	&=f_{0}^{(n)}\wedge \frac{\partial}{\partial z}f_{0}^{(n)}
    	\wedge\cdots\wedge \frac{\partial^q}{\partial z^q}f_{0}^{(n)}\\
    	&=f_{0}^{(n)}\wedge f_{1}^{(n)}\wedge\cdots\wedge f_{q}^{(n)}
    \end{aligned}
    \end{equation}
    represent the \(q\)-th osculating Pl\"ucker curve. At a
    ramification point, cancel the common zeros of the components
    of \(F_q^{(n)}\) and write \(\widehat F_q^{(n)}\) for the resulting
    nowhere-zero local holomorphic lift. Its degree is
    \(\mathfrak d_q^{(n)}\), where
    \begin{equation}\label{eq:osculating-degree-integral}
    	\begin{aligned}
    	&\mathfrak d_q^{(n)}
    	=\frac{1}{2\pi\mathrm{i}}\int_{\mathbb{S}^2}
    	\partial\bar\partial\log |\widehat F_q^{(n)}|^2
    	\mathrm{d}\overline{z} \wedge \mathrm{d}z\\
    	&\qquad =\frac{1}{2\pi\mathrm{i}}
    	\int_{\mathbb{S}^2}\lambda_q^{(n)}
    	\,\mathrm{d}\overline{z} \wedge \mathrm{d}z.
    \end{aligned}
    \end{equation}
    Here \(0\leq q<n\); the top osculating curve is constant,
    so \(\mathfrak d_n^{(n)}=0\).
    The first integrand is the globally defined curvature density
    computed from these local saturated lifts; the squared norm
    density on the right extends smoothly as the squared norm
    of the actual Frenet connection. In the unramified case no
    cancellation is necessary. Integrating the logarithm of an
    unsaturated lift distributionally at a common zero would instead
    include its zero-divisor contribution.

    For the seed in the replacement construction below, choose a
    primitive polynomial lift \(f_0\) of arbitrary degree and
    suppress the superscript in its monic derivative lifts
    \(f_0,\ldots,f_m\). Their forward and backward equations are those
    of Lemma~\ref{lem:veronese-frenet}, with the conjugate equations
    obtained by complex conjugation. We suppress the remaining
    superscript and write \(\lambda_i=\lambda_i^{(m)}\). We also use
    \[
     \mathscr F_p=\bigoplus_{i=0}^p\mathcal L_i,\qquad
     \mathscr F_{-1}=0,
    \]
    with \(\lambda_{-1}=\lambda_m=0\). Here \(\mathscr F_p\)
    is a prefix bundle, whereas \(F_p\) denotes the exterior
    vector in \eqref{eq:osculating-curve}. Each \(\lambda_i\)
    is a smooth nonnegative connection density. If \(\epsilon_i\)
    is the degree of the zero divisor of the Frenet homomorphism
    \(\mathcal L_i\to\mathcal L_{i+1}\otimes K_{\mathbb S^2}\),
    then \(\deg K_{\mathbb S^2}=-2\) gives
    \begin{equation}\label{eq:ramified-frenet-degrees}
     \deg\mathcal L_{i+1}=\deg\mathcal L_i+2+\epsilon_i.
    \end{equation}

    If the entire Frenet sequence of \(f_0^{(n)}\) is totally
    unramified, then
    \begin{equation}\label{eq:unramified-degree-formula}
    	\mathfrak d_q^{(n)}=(q+1)(n-q).
    \end{equation}
    Let
    \(V_q^{(n)}=(v_{q,0},v_{q,1},\dots,v_{q,n})^{\transp}\), where
    \begin{equation}
    	v_{q,r}(z)
    	=\frac{q!}{(1+z\overline{z})^q}
    	\sqrt{\binom{n}{r}}z^{r-q}
    	\sum_{k=0}^{q}
    	(-1)^k\binom{r}{q-k}\binom{n-r}{k}
    	(z\overline{z})^k,
    \end{equation}
    The vectors \(V_q^{(n)}\) form the Veronese harmonic sequence:
    
    \begin{equation}
    	0\stackrel{\partial''}{\leftarrow}\underline{V_0^{(n)}}
    	\stackrel{\partial'}{\rightarrow}\cdots
    	\stackrel{\partial'}{\rightarrow}\underline{V_q^{(n)}}
    	\stackrel{\partial'}{\rightarrow}\cdots
    	\stackrel{\partial'}{\rightarrow}\underline{V_n^{(n)}}
    	\stackrel{\partial'}{\rightarrow}0.
    \end{equation}
    Here
    \(\lVert V_q^{(n)}\rVert^2
    =\frac{n!q!}{(n-q)!}(1+z\overline{z})^{n-2q}\), and the metric
    induced by \(V_q^{(n)}\) is
    \begin{equation}
    	\mathrm{ds}_q^2=\frac{n+2q(n-q)}{(1+z\overline{z})^2}\mathrm{d}z\mathrm{d}\overline{z},
    \end{equation}
    with curvature
    \begin{equation}
    	K_q=\frac{4}{n+2q(n-q)}.
    \end{equation}
    For comparison, we recall the projective rigidity theorem of
    Bolton, Jensen, Rigoli, and Woodward.
    
    \begin{theorem}[Veronese rigidity in projective space]
    	\label{thm:veronese-rigidity}
    	\citep{bolton1988conformal}
    	Suppose that \(\varphi:\mathbb{S}^2\to\mathbb{CP}^n\) is a linearly
    	full conformal minimal immersion of constant Gaussian curvature.
    	Then there is a spherical affine coordinate in which, up to
    	a constant unitary transformation of
    	\(\mathbb{CP}^n\), the harmonic sequence induced by \(\varphi\) is
    	the Veronese sequence.
    \end{theorem}

    \section{Replacement Towers and Short-Window Graphs}\label{sec:basic-setup}
    Every irreducible harmonic map
    \(\varphi:\mathbb{S}^2\to G(2,n;\mathbb{R})\) of finite isotropy order
    \(r=2s+1\) admits the construction in
    Lemma~\ref{lem:bew-replacement-series}.
    \begin{lemma}[Replacement tower]
    	\label{lem:bew-replacement-series}
    	\citep{bahy1989explicit}
    	There are unique harmonic maps
    	\(V_j:\mathbb{S}^2\to G(2,n;\mathbb{R})\)
    	(\(j=0,1,\dots,\ell\)) and
    	\(U_j:\mathbb{S}^2\to G(2,n;\mathbb{C})\)
    	(\(j=0,1,\dots,\ell-1\)) such that:
    	
    	(1) \(V_0\) is a real mixed pair: \(f_0^{(m)}\in H_n^{\ell+s}\) and
    	\(\underline{V_0}
    	=\underline{f_0^{(m)}}\oplus
    	\underline{\overline{f_0}^{(m)}}\).
    	
    	(2) For \(j=0,1,\dots,\ell-1\), the map \(U_j\) is obtained from
    	\(V_j\) by forward replacement along \(\underline{\rho_j}\); that
    	is,
    	\begin{equation}\underline{U_{j}}=(\underline{\rho_j}^\perp \cap \underline{V_{j}}) \oplus \operatorname{Im}(A'_{V_j}|\underline{\rho_j}),
    	\end{equation}
    	where \(\underline{\rho_j}\) is the image of the
    	\(\partial''\)-closing map of \(V_j\).
    	
    	(3) For \(j=0,1,\dots,\ell-1\), the map \(V_{j+1}\) is obtained from
    	\(U_j\) by backward replacement along
    	\(\underline{\gamma_j}^{\perp}\cap\underline{U_j}\). Thus
    	\begin{equation}
    		\underline{V_{j+1}}
    	=\underline{\gamma_j}\oplus
    	\operatorname{Im}\left(A''_{U_{j}}|\underline{\gamma_j}^{\perp}\cap \underline{U_{j}}\right),
    	\end{equation}
    	where \(\underline{\gamma_j}\subset U_j\) is a holomorphic subbundle
    	that is not the image of the \(\partial'\)-closing map of \(U_j\).
    	
    	(4) \(V_\ell=\varphi\).
    	
    	(5) For every \(j=0,1,\dots,\ell\), the map
    	\(V_j:\mathbb{S}^2\to G(2,n;\mathbb{R})\) is harmonic.
    	
    	(6) For \(j=0,1,\dots,\ell-1\), the map \(U_j\) is harmonic and
    	satisfies
    	\(\partial''(\underline{U_j})
    	=\underline{\overline{U_j}}\); consequently,
    	\(\underline{V_{j+1}}
    	=\underline{\gamma_j}\oplus
    	\underline{\overline{\gamma_j}}\).
    	
    	(7) For \(j=0,1,\dots,\ell\), the isotropy order of \(V_j\) is
    	\(2(\ell+s-j)+1\). For \(j=0,1,\dots,\ell-1\), the isotropy order of
    	\(U_j\) is \(2(\ell+s-j)\).
    	
    	The class \(H_n^{\ell+s}\) consists of all linearly full holomorphic maps
    	\(f_0^{(m)}:\mathbb{S}^2\to\mathbb{CP}^m
    	\hookrightarrow\mathbb{CP}^{n-1}\), whose image lies in
    	\(Q_{n-2}\) and satisfies
    	\begin{equation}
    		\begin{cases}
    		b(f_0^{(m)},f_i^{(m)})=0,
    		\quad 0\leq i\leq2(\ell+s)+1,  \\
    		b(f_0^{(m)},f_{2(\ell+s)+2}^{(m)})\not\equiv0.
    	\end{cases}.
    	\end{equation}
    	\end{lemma}

    The other isotropic lift is its complex conjugate. It produces the
    conjugate replacement tower, the same real layers, and the same
    actual depth.
    
    Throughout, \(m\) denotes the dimension of \(\mathbb{CP}^m\),
    \(n\) is the dimension parameter of \(G(2,n;\mathbb{R})\), \(\ell\) is the
    number of forward-backward replacements, and \(r=2s+1\) is the isotropy
    order of \(\varphi\). The maps \(V_j,U_j,\gamma_j,\rho_j\) are the
    intermediate objects in the replacement construction. We abbreviate
    \(f_0^{(m)},\dots,f_m^{(m)}\) by \(f_0,\dots,f_m\), and we identify
    each map with its induced bundle. For a complex line subbundle
    \(\underline L\), the symbol \(L\) without an underline denotes a
    chosen local lift on the regular set. When a formula depends on that
    choice, we state the normalization explicitly. Terms with indices
    outside their displayed ranges are zero.

    \begin{definition}[Real mixed pairs and real layers]
    	\label{def:real-mixed-layer}
    	\citep{bahy1989explicit}
    	Let \(f:M\to Q_{n-2}\) be a holomorphic map, and let
    	\(\mathcal F\) denote its associated holomorphic line subbundle of
    	\(\underline{\mathbb{C}^n}\). A \emph{real mixed pair} is the real
    	two-plane bundle
    	\[
    		\mathcal F\oplus\overline{\mathcal F}.
    	\]
    	Thus the two summands of a real mixed pair are a holomorphic
    	subbundle and its complex conjugate. In contrast, for an arbitrary
    	complex line \(L\subset Q_{n-2}\), we call
    	\(L\oplus\overline L\) its \emph{real layer}; the summands need not
    	be holomorphic. Every real mixed pair is a real layer, but the
    	converse is false.

    	The maps \(V_j\) in Lemma~\ref{lem:bew-replacement-series} are
    	\emph{real layers} of the replacement tower. By part (1) of that
    	lemma, \(V_0\) is the initial real mixed pair; no such claim is
    	made for \(V_j\) when \(j>0\). The target map is the \(\ell\)-th
    	real layer \(V_\ell=\varphi\).
    \end{definition}

    \begin{corollary}[BEW admissible range]
    	\label{cor:bew-admissible-range}
    	Under the hypotheses of Lemma~\ref{lem:bew-replacement-series},
    	\begin{equation}\label{eq:bew-admissible-range}
    		2(\ell+s)+2\leq m.
    	\end{equation}
    \end{corollary}

    \begin{proof}
    	The class \(H_n^{\ell+s}\) requires the nonzero Frenet vector
    	\(f_{2(\ell+s)+2}^{(m)}\). Since the Frenet sequence of the curve in
    	\(\mathbb{CP}^m\) is indexed by \(0,\ldots,m\), this vector cannot
    	occur unless \(2(\ell+s)+2\leq m\).
    \end{proof}

Differentiating the defining zero pairings of
\(H_n^{\ell+s}\) gives
\begin{equation}\label{eq:short-mixed-orthogonality}
 b(f_i,f_h)=0\quad(i+h\leq2(\ell+s)+1=2\ell+r).
\end{equation}
The same equations hold for the saturated lines. Thus
\(\mathscr F_p\perp_H\overline{\mathscr F_q}\) when
\(p+q\leq2\ell+r\). Thus every graph calculation below uses an
orthogonal short window; no assertion is made for a larger fixed window.

\begin{table}[H]
\centering\small
\caption{Notation retained throughout the argument.}
\label{tab:main-notation}
\begin{tabular}{@{}lp{0.73\textwidth}@{}}
\toprule
Symbol & Meaning \\
\midrule
\(r=2s+1,\ \ell+s\) & target isotropy order and seed class index \\
\(V_j,U_j\) & real even layer and complex odd layer \\
\(\rho_j,\gamma_j\) & closing line and retained holomorphic line in the replacement construction \\
\(v_j,u_j,x_j\) & retained lift, forward vector, and deviation \\
\(\mathcal L_i,\mathscr F_p,F_j\) & saturated Frenet line, prefix bundle, and exterior vector \\
\(\lambda_i,\lambda_i^{(n)}\) & projective Frenet connection density \\
\(S^j,w_i^j,\mathbb L_j\) & graph matrix, graph basis, and graph subspace \\
\(\Psi_{j,\nu},\boldsymbol\Psi_j,P_j\) & gradings, raw augmented vector, and raw norm \\
\(\widehat{\boldsymbol\Psi}_j,Q_j,D_j\) & independently saturated vector, norm, and degree \\
\(\mathsf R,\mathcal G(\mathsf R)\) & complex symmetric Gram matrix and real Gram matrix \\
\bottomrule
\end{tabular}
\end{table}

    The actual depth is read from the noncancelling construction of
\citep[Theorem~4.7 and \S5(A)]{bahy1989explicit}. The retained line at
level \(j\geq1\) is the line \(\gamma_{j-1}\) of
Lemma~\ref{lem:bew-replacement-series}; it has a nonzero new highest
Frenet component generically. Normalizing that component to one is
legitimate on a dense set and does not assert that it has no zeros
globally. Such zeros will be measured by a divisor in
Section~\ref{sec:degree-growth}.

For completeness, the generic highest-row assertion follows
inductively from the two half-steps. If the retained lift at level
\(j\) has highest positive row \(f_j\), its forward normal
derivative has highest row \(f_{j+1}\) with the same nonzero
coefficient: the current projection has no row of that index.
The other line of the odd layer is \(\overline{\underline v_j}\),
which has no such positive row. The decreasing-replacement
construction identifies that unchanged line with the excluded
line of the odd layer. A retained line other than
it therefore has a nonzero forward component generically and
can be normalized to \(u_j+x_j\overline{v_j}\). Starting at
\(v_0=f_0\) proves the assertion at every level. This uses
noncancellation of the canonical decomposition, not global
nonvanishing of a derivative or of a projection.

\noindent\textbf{External inputs.}
For orientation, the construction uses external results at three
separate points. Theorem~4.7 of \citep{bahy1989explicit} supplies
existence and uniqueness of the noncancelling replacement tower.
The replacement formulas and induced holomorphic structures give the
actual forward and backward line maps. Proposition~5.1, Lemma~5.2,
and Theorem~5.3 supply the rational representatives used later.
All graph, Pfaffian, degree, and valuation calculations after these
inputs are internal to the present paper. In particular, no
unramifiedness or constant-curvature classification is hidden in the
cited construction.

\subsection{Retained lines and deviations}\label{sec:deviation-structure}

We now translate the two geometric half-steps into scalar data.
All monic formulas in this and the next section are first read on
the common regular set; no conclusion about ramification is built
into the notation.

    	At the \(j\)-th stage, let \(\underline v_j\) denote the retained
    	line subbundle. On a graph-regular chart choose a local
    	nowhere-zero lift \(v_j\), normalized so that its coefficient
    	along \(f_j\) is \(1\) whenever this coefficient is nonzero. All
    	formulas in this section use this chosen lift. The quantities
    	\(p_i^j,q_i^j,x_j,S^j\) below are therefore defined relative to
    	this lift; for statements about the underlying line, we check
    	explicitly whether they depend on this choice. Set
    	\begin{equation}\label{eq:deviation-forward}
    		u_j
    		=\partial v_j
    		-\bigl(\partial\log\lVert v_j\rVert^2\bigr)v_j.
    	\end{equation}
    	For this monic compatible frame, the graph connection and
    	orthogonality imply \(H(\bar\partial v_j,v_j)=0\); this is
    	proved in Appendix~\ref{app:graph-connection}. Consequently,
    	\eqref{eq:deviation-forward} is precisely the normal projection
    	of \(\partial v_j\). For an arbitrary smooth frame one should
     use the orthogonal projection onto the normal bundle instead.
    	The next retained line and the associated layers are
    	\begin{align}
    		v_{j+1} &= u_j + x_j\overline{v_j},
    		\label{eq:deviation-replacement}\\
    		V_j&=v_j\oplus\overline{v_j},\qquad
    		U_j=u_j\oplus\overline{v_j},
    		\label{eq:deviation-real-bundles}
    	\end{align}
    	where \(x_j\) is a smooth complex-valued function on the chosen
    	chart. Its value depends on the coordinate and the lift, but its
    	vanishing is invariant. Here \(0\leq j<\ell\); the formula for \(V_j\) also
    	applies at \(j=\ell\). Equivalently,
    	\begin{equation}\label{eq:paper-even-recurrence}
    		v_j=u_{j-1}+x_{j-1}\overline{v_{j-1}},
    		\qquad 1\leq j\leq \ell.
    	\end{equation}
    
    \begin{definition}\label{def:deviation}
    	We call \(x_j\) the deviation parameter. For \(0\leq j<\ell\), the
    	forward-backward transformation from \(v_j\) to \(v_{j+1}\), with
    	parameter \(x_j\), is faithful when \(x_j=0\). If this holds for
    	every \(0\leq j<\ell\), the replacement tower ending at \(V_\ell\) is
    	totally faithful.
    \end{definition}
    
    \begin{proposition}[Faithful towers and Frenet lines]
    	\label{prop:faithful-frenet}
    	The replacement tower is totally faithful if and only if its
    	normalized lifts satisfy \(v_j=f_j\) for
    	\(0\leq j\leq\ell\). In this case, the retained complex lines
    	form the first \(\ell+1\) members of the projective harmonic
    	sequence generated by \(f_0\).
    \end{proposition}

    \begin{proof}
    	Starting with \(v_0=f_0\), if \(x_j=0\), then
    	\eqref{eq:deviation-forward} and the normalized Frenet equation
    	give \(v_{j+1}=f_{j+1}\) whenever \(v_j=f_j\). Induction proves
    	the forward implication. Conversely, if \(v_j=f_j\) for every
    	\(j\), then \eqref{eq:deviation-replacement} reads
    	\(f_{j+1}=f_{j+1}+x_j\overline{f_j}\). Since \(f_j\) is nowhere
    	zero, \(x_j=0\). The last assertion is the Frenet construction
    	for the initial holomorphic curve in its projective span
    	\(\mathbb{CP}^m\).
    \end{proof}

    \subsection{A New Coordinate System}
    
    We now introduce a more convenient coordinate system. Write
    \begin{equation}\label{eq:v-expansion}
    	v_j=\sum_{i=0}^{j}p_i^j f_i+
    	\sum_{i=0}^{j-1}q_i^j\overline{f_i},
    	\qquad p_j^j=1,\quad q_j^j=0.
    \end{equation}
    The normalization \(p_j^j=1\) follows from
    \eqref{eq:deviation-forward} and \eqref{eq:deviation-replacement};
    indeed, the coefficient of the highest Frenet vector \(f_j\) equals
    one.
    
    \begin{lemma}[Coefficient recurrence for replacement layers]
    	\label{lem:coefficient-recurrence}
    	The coefficients in \eqref{eq:v-expansion} satisfy
    	\begin{align}
    		p_i^{j+1}
    		&= \partial p_i^j+
    		\left(\partial\log\norm{f_i}^2-\partial\log\norm{v_j}^2\right)p_i^j
    		+p_{i-1}^j+x_j\overline{q_i^j},\label{eq:p-recurrence}\\
    		q_i^{j+1}
    		&= \partial q_i^j-\left(\partial\log\norm{v_j}^2\right)q_i^j
    		-\lambda_iq_{i+1}^j
    		+x_j\overline{p_i^j}.\label{eq:q-recurrence}
    	\end{align}
    \end{lemma}
    
    \begin{proof}
    	Equations~\eqref{eq:p-recurrence} and~\eqref{eq:q-recurrence} follow
    	by substituting the expansion of \(v_j\) into the recurrence
    	\begin{equation}
    		v_{j+1}= \partial v_j-
    	\left(\partial\log\norm{v_j}^2\right)v_j+x_j\overline{v_j}
    	\end{equation}
    	and comparing the coefficients of \(f_i\) and
    	\(\overline{f_i}\).
    	For example, collect the coefficient of \(f_i\). Differentiating
    	\(p_i^j f_i\) gives
    	\[
    		\partial(p_i^j f_i)
    		=\partial p_i^j\,f_i
    		+p_i^jf_{i+1}
    		+p_i^j\bigl(\partial\log\norm{f_i}^2\bigr)f_i.
    	\]
    	The \(f_i\)-component of
    	\(\partial(p_{i-1}^j f_{i-1})\) is \(p_{i-1}^j f_i\), while the
    	tangential subtraction contributes
    	\(-(\partial\log\norm{v_j}^2)p_i^j f_i\). Finally,
    	\(x_j\overline{v_j}\) contributes \(x_j\overline{q_i^j}f_i\).
    	Adding these terms proves \eqref{eq:p-recurrence}; absent indices
    	are omitted at the endpoints.

    	For the coefficient of \(\overline{f_i}\), use
    	\(\partial\overline{f_{i+1}}
    	=-\lambda_i\overline{f_i}\). The derivative
    	\(\partial(q_i^j\overline{f_i})\) therefore contributes
    	\(\partial q_i^j\) at \(\overline{f_i}\), while the
    	\(i+1\)-term contributes \(-\lambda_iq_{i+1}^j\). The remaining
    	contributions are
    	\(-(\partial\log\norm{v_j}^2)q_i^j\) from the tangential
    	subtraction and \(x_j\overline{p_i^j}\) from the replacement.
    	This proves \eqref{eq:q-recurrence}.
    \end{proof}
    
    The vectors \(v_0,\dots,v_j\) span a \((j+1)\)-dimensional space,
    whereas the ambient Frenet vectors
    \begin{equation}
    	f_0,\ldots,f_j,\overline{f_0},\ldots,\overline{f_{j-1}}
    \end{equation}
    span a space of dimension \(2j+1\). To eliminate the conjugate
    components in \eqref{eq:v-expansion}, we construct a basis
    \(w_0^j,\dots,w_j^j\) that absorbs them while preserving the information
    carried by \(v_j\), so that
    \begin{equation}\label{eq:v-in-w}
    	v_j=\sum_{i=0}^{j}p_i^j w_i^j.
    \end{equation}
    We first define the skew-symmetric matrices that determine this basis.
    
    \begin{definition}\label{def:S-and-w}
    	Define skew-symmetric matrices \(S^j\) inductively by
    	\begin{equation}
    		S^0=(0),\qquad
    	S^j=
    	\begin{pmatrix}
    		S^{j-1}&t^j\\
    		-(t^j)^{\transp}&0
    	\end{pmatrix},
    	\qquad
    	t^j=
    	\begin{pmatrix}
    		t_0^j\\ \vdots\\ t_{j-1}^j
    	\end{pmatrix},
    	\end{equation}
    	where
    	\begin{equation}\label{eq:t-definition}
    		t_i^j=\norm{f_i}^2q_i^j+
    		\sum_{k=0}^{j-1}p_k^j S_{ki}^{j-1},
    		\qquad 0\leq i<j.
    	\end{equation}
    	For \(0\leq i\leq j\), set
    	\begin{equation}\label{eq:w-definition}
    		w_i^j=f_i+\sum_{k=0}^{j}
    		\frac{S_{ki}^j}{\norm{f_k}^2}\overline{f_k}.
    	\end{equation}
    \end{definition}
    
    \begin{note}
    	The upper-left block of \(S^j\) is \(S^{j-1}\). When the level is
    	clear, we omit the superscript and write \(S_{ki}\) instead of
    	\(S_{ki}^j\). The superscript on \(w_i^j\), however, must be retained,
    	since in general \(w_i^{j+1}\ne w_i^j\). In fact,
    	\begin{equation}\label{eq:w-update}
    		w_i^j=w_i^{j-1}
    		-\frac{t_i^j}{\norm{f_j}^2}\overline{f_j},
    		\qquad 0\leq i<j.
    	\end{equation}
    	Consequently, vectors from different levels should not be used
    	together in one basis expansion.
    \end{note}

    \begin{example}\label{ex:graph-l1}
    	For \(\ell=1\), the coefficient recurrence gives
    	\begin{equation}
    		p_0^1=0,\qquad q_0^1=x_0.
    	\end{equation}
    	Since \(S^0=(0)\), Definition~\ref{def:S-and-w} gives
    	\begin{equation}
    		t_0^1=\norm{f_0}^2q_0^1
    		=x_0\norm{f_0}^2.
    	\end{equation}
    	Hence the complete matrix and graph basis are
    	\begin{equation}
    		S^1=
    		\begin{pmatrix}
    			0&x_0\norm{f_0}^2\\
    			-x_0\norm{f_0}^2&0
    		\end{pmatrix},
    		\qquad
    		\begin{pmatrix}
    			w_0^1\\ w_1^1
    		\end{pmatrix}
    		=
    		\begin{pmatrix}
    			f_0-\dfrac{x_0\norm{f_0}^2}{\norm{f_1}^2}
    			\overline{f_1}\\[2mm]
    			f_1+x_0\overline{f_0}
    		\end{pmatrix}.
    	\end{equation}
    	Following \citep{jiao2017classification,li2023conformal}, we note
    	that the quantity \(x_0\norm{f_0}^2\), also occurring in
    	Examples~\ref{ex:graph-l2}, \ref{ex:pfaffian-l1}, and
    	\ref{ex:pfaffian-l2}, is denoted by \(h\) in these references.
    \end{example}

    \begin{example}\label{ex:graph-l2}
    	For \(\ell=2\), substituting the coefficient recurrence into the
    	definition of \(t_i^2\) gives
    	\begin{align}
    		t_0^2
    		&=\partial\bigl(x_0\norm{f_0}^2\bigr)
    		-x_0\norm{f_0}^2
    		\partial\log\left(
    		\norm{f_0}^2\norm{f_1}^2\right)=:A,\\
    		t_1^2
    		&=x_1\left(
    		\norm{f_1}^2+|x_0|^2\norm{f_0}^2\right)
    		=:B.
    	\end{align}
    	Thus the complete level-two matrix is
    	\begin{equation}\label{eq:example-S2}
    		\boxed{
    		S^2=
    		\begin{pmatrix}
    			0&x_0\norm{f_0}^2&A\\
    			-x_0\norm{f_0}^2&0&B\\
    			-A&-B&0
    		\end{pmatrix}.}
    	\end{equation}
    	The corresponding graph basis is
    	\begin{equation}\label{eq:example-w2}
    		\begin{aligned}
    			w_0^2&=f_0-
    			\frac{x_0\norm{f_0}^2}{\norm{f_1}^2}\overline{f_1}
    			-\frac{A}{\norm{f_2}^2}
    			\overline{f_2},\\
    			w_1^2&=f_1+x_0\overline{f_0}
    			-\frac{B}{\norm{f_2}^2}\overline{f_2},\\
    			w_2^2&=f_2+
    			\frac{A}{\norm{f_0}^2}
    			\overline{f_0}
    			+\frac{B}{\norm{f_1}^2}\overline{f_1}.
    		\end{aligned}
    	\end{equation}
    \end{example}
    
    \begin{lemma}[Graph decomposition and pairing]
    	\label{lem:graph-properties}
    	Let \(\mathbb L_j=\operatorname{span} \{w_0^j,\ldots,w_j^j\}\). Then
    	\begin{equation}\label{eq:graph-decomposition}
    		\operatorname{span}\{f_0,\ldots,f_j,\overline{f_0},\ldots,\overline{f_j}\}
    		=\mathbb L_j\oplus_H\overline{\mathbb L_j},
    	\end{equation}
    	where \(\oplus_H\) denotes an orthogonal direct sum with respect to
    	the Hermitian form \(H\). Moreover,
    	\begin{equation}\label{eq:t-pairing}
    		t_i^j=b(v_j,w_i^{j-1}),\qquad 0\leq i<j,
    	\end{equation}
    	and the representation \eqref{eq:v-in-w} holds.
    \end{lemma}
    
    \begin{proof}
    	We first show that \(\mathbb L_j\) is \(b\)-isotropic. By
    	\eqref{eq:w-definition} and the Frenet pairing relations,
    	\begin{align}
    		b(w_i^j,w_k^j)
    		&=S_{ik}^j+S_{ki}^j\\
    		&=0.
    	\end{align}
    	Since \(H(w_i^j,\overline{w_k^j})=b(w_i^j,w_k^j)\), we obtain
    	\(\mathbb L_j\perp_H\overline{\mathbb L_j}\).
    	
    	The non-conjugate component of \(w_i^j\) is exactly \(f_i\), so
    	\(w_0^j,\ldots,w_j^j\) are linearly independent. Their conjugates
    	are also linearly independent. Thus both summands have dimension
    	\(j+1\), and their Hermitian orthogonality implies that their
    	intersection is zero. Their direct sum consequently has dimension
    	\(2j+2\), which is the dimension of the space on the left-hand side
    	of \eqref{eq:graph-decomposition}. This proves
    	\eqref{eq:graph-decomposition}.
    	
    	It remains to prove the representation of \(v_j\). By
    	\eqref{eq:t-definition},
    	\begin{equation}
    		t_k^j=\norm{f_k}^2q_k^j+
    	\sum_{r=0}^{j-1}p_r^jS_{rk}^{j-1}.
    	\end{equation}
    	Using skew-symmetry, this is equivalent to
    	\begin{equation}\label{eq:q-k-inversion}
    		\frac{t_k^j+\sum_{r=0}^{j-1}p_r^jS_{kr}^{j-1}}
    	{\norm{f_k}^2}=q_k^j.
    	\end{equation}
    	For \(k<j\), the coefficient of \(\overline{f_k}\) in
    	\begin{equation}
    		w_j^j+\sum_{i=0}^{j-1}p_i^j w_i^j
    	\end{equation}
    	is precisely the left-hand side of \eqref{eq:q-k-inversion}, whereas
    	the coefficient of \(f_k\) is \(p_k^j\). The coefficient of \(f_j\) is
    	one, and the coefficient of \(\overline{f_j}\) is
    	\begin{equation}
    		\frac{1}{\norm{f_j}^2}\sum_{i=0}^{j-1}p_i^jS_{ji}^j
    	=-\frac{1}{\norm{f_j}^2}\sum_{i=0}^{j-1}p_i^jt_i^j.
    	\end{equation}
    	This final sum vanishes because
    	\begin{equation}
    		\sum_{i=0}^{j-1}p_i^jt_i^j
    	=\sum_{i=0}^{j-1}\norm{f_i}^2p_i^jq_i^j
    	+\sum_{i,k=0}^{j-1}p_i^jp_k^jS_{ki}^{j-1}=0.
    	\end{equation}
    	The first term vanishes by the isotropy condition
    	\(b(v_j,v_j)=0\) and \(q_j^j=0\); the second vanishes because
    	\(S^{j-1}\) is skew-symmetric. Hence
    	\(v_j=\sum_{i=0}^jp_i^jw_i^j\), as required. Finally,
    	\eqref{eq:t-pairing} follows directly from \eqref{eq:t-definition}
    	and the definition of \(w_i^{j-1}\) in
    	\eqref{eq:w-definition}.
    \end{proof}
    
    We first record a differential consequence of the graph connection
    equation.
    
    \begin{lemma}[Backward derivative under graph connection]
    	\label{lem:conditional-graph-connection}
    Suppose that \(S^j\) satisfies
    \begin{equation}\label{eq:S-connection-assumption}
    	 \bar\partial S_{ki}^j
    	=-\lambda_{k-1}S_{k-1,i}^j-\lambda_{i-1}S_{k,i-1}^j,
    \end{equation}
    Then
    \begin{equation}\label{eq:w-backward-derivative}
    		 \bar\partial w_i^j
    		=-\lambda_{i-1}w_{i-1}^j
    		+\frac{S_{ji}^j}{\norm{f_j}^2}\overline{f_{j+1}}.
    \end{equation}
    \end{lemma}
    
    \begin{proof}
    	By the Frenet equations, the derivative of each term in
    	\eqref{eq:w-definition} is
    	\begin{equation}
    		\bar\partial\left(
    	\frac{S_{ki}^j}{\norm{f_k}^2}\overline{f_k}
    	\right)
    	=
    	\frac{\bar\partial S_{ki}^j}{\norm{f_k}^2}\overline{f_k}
    	+\frac{S_{ki}^j}{\norm{f_k}^2}\overline{f_{k+1}}.
    	\end{equation}
    	Indeed, the derivative of \(1/\norm{f_k}^2\) cancels the
    	\(\log\norm{f_k}^2\) term in \(\overline{f_k}\). Substituting
    	\eqref{eq:S-connection-assumption}, the first contribution is
    	\begin{equation}
    		-\sum_{k=1}^{j}
    	\frac{S_{k-1,i}^j}{\norm{f_{k-1}}^2}\overline{f_k}
    	-\sum_{k=0}^{j}
    	\lambda_{i-1}\frac{S_{k,i-1}^j}{\norm{f_k}^2}\overline{f_k}.
    	\end{equation}
    	The first sum cancels the terms
    	\(\overline{f_1},\dots,\overline{f_j}\) produced by the second part
    	of the product rule. Only the top term
    	\begin{equation}
    		\frac{S_{ji}^j}{\norm{f_j}^2}\overline{f_{j+1}}
    	\end{equation}
    	remains. Because
    	\begin{equation}
    		\bar\partial f_i
    	=-\lambda_{i-1}f_{i-1},
    	\end{equation}
    	the uncancelled lower-index terms combine to
    	\begin{equation}
    		-\lambda_{i-1}w_{i-1}^j.
    	\end{equation}
    	This proves \eqref{eq:w-backward-derivative}.
    \end{proof}
    
    The graph matrices arising from the replacement sequence satisfy
    compatible orthogonality, connection, and contact identities.
    We state these together with the auxiliary derivative identities
    needed in \S\ref{sec:rigidity}; their simultaneous-induction proof
    is given in Appendix~\ref{app:graph-connection}.
    
    For clarity, set
    \begin{equation}
    	C_{ik}^j=b(\partial w_i^j,w_k^j).
    \end{equation}
    The three basic graph identities are as follows.
    
    \begin{lemma}[Graph connection and contact identities]\label{lem:simultaneous}
    	For every \(0\leq j\leq \ell\), the following statements hold:
    	
    		(1) \(v_j\in\mathbb L_j\), and
    		\begin{equation}\label{eq:orthogonality-v-w}
    			H(v_j,w_i^j)=\norm{v_j}^2\delta_{ij},
    			\qquad 0\leq i\leq j.
    		\end{equation}
    		(2) The matrix \(S^j\) satisfies
    		\begin{equation}\label{eq:S-connection}
    			 \bar\partial S_{ki}^j
    			=-\lambda_{k-1}S_{k-1,i}^j
    			-\lambda_{i-1}S_{k,i-1}^j,
    			\qquad 0\leq i,k\leq j,
    		\end{equation}
    		(3) The contact coefficients satisfy
    		\begin{equation}\label{eq:contact-condition}
    			C_{ik}^j=0,\qquad 0\leq i,k<j.
    		\end{equation}

        For \(1\leq j\leq\ell\), the new column is also given by
        \begin{equation}\label{eq:contact-new-column}
         t_i^j=C_{j-1,i}^{j-1}
         =b(\partial w_{j-1}^{j-1},w_i^{j-1}),
         \qquad0\leq i<j-1.
        \end{equation}

        In addition, use \eqref{eq:deviation-forward} to define \(u_j\)
        also at \(j=\ell\). For \(0\leq j<\ell\), there exist
        \(Z_j\in\mathbb L_j\) and a scalar function \(c_j\) such that
        \begin{equation}\label{eq:Ubar-decomposition}
            \bar\partial u_j=Z_j+c_j\overline{v_j},
        \end{equation}
        and the deviation parameter satisfies
        \begin{equation}\label{eq:scalar-replacement}
            c_j+\bar\partial x_j
            +x_j\bar\partial\log\norm{v_j}^2=0.
        \end{equation}
        At the terminal level,
        \begin{equation}\label{eq:terminal-graph-connection}
            \begin{aligned}
            \bar\partial v_\ell&\in
            \operatorname{span}\{w_0^\ell,\ldots,w_{\ell-1}^\ell\},\\
            \bar\partial u_\ell&=Z_\ell+c_\ell\overline{v_\ell},
            \quad Z_\ell\in\mathbb L_\ell.
            \end{aligned}
        \end{equation}
        Moreover, \(u_\ell\perp_H\mathbb L_\ell\) and
        \(b(u_\ell,v_\ell)=0\).
    \end{lemma}
    
    \begin{proof}[Proof outline]
        The old connection equation gives the backward derivative of the
        retained lift, and the old contact condition restricts the
        derivative of its forward direction. Induced holomorphicity then
        gives \eqref{eq:scalar-replacement}. The graph update yields
        orthogonality at the new level, followed by its connection
        equation and contact condition, in that order. The terminal
        identities use only the existing level-\(\ell\) graph, with no
        further replacement. The full simultaneous-induction argument
        is given in Appendix~\ref{app:graph-connection}.
    \end{proof}

\begin{note}[Regular set and terminal level]
The monic graph identities hold on the common regular set. Across
excluded points we use the smooth saturated bundles and actual
connection energies. The induction in Appendix~\ref{app:graph-connection} uses
mixed index sums at most \(2\ell+1\), within
\eqref{eq:short-mixed-orthogonality}, including when \(s=0\).
At the terminal level its inputs are
\(u_\ell\perp_H\mathbb L_\ell\) and \(b(u_\ell,v_\ell)=0\).
\end{note}

\section{Rational Pfaffian Data and Independent Saturation}
\label{sec:pfaffian-data}

The graph equations give holomorphic Pfaffian vectors whose curvature
forms measure the actual replacement energies. We derive their norm,
connection, and source identities in the monic gauge, then saturate
each rational curve independently. The raw quantities are \(P_j,B_j\);
the saturated quantities are \(Q_j,D_j\) and the hatted vectors.

    \subsection{Pfaffian}\label{subsec:pfaffian}

    We use elementary Pfaffian identities, with proofs included to
    fix the signs. Their role here is to package the graph connection
    and contact conditions into uniformly holomorphic graded data and
    exact norm and mixed-kernel formulas at every replacement depth.

    \begin{definition}\label{def:pfaffian-gradings}
    	For a \(2m\times2m\) complex skew-symmetric matrix
    	\(A=(a_{ij})\), define its Pfaffian by
    	\begin{equation}\label{eq:matrix-pfaffian}
    		\operatorname{Pf}(A)
    		=
    		\frac{1}{2^m m!}
    		\sum_{\sigma\in\mathfrak S_{2m}}
    		\operatorname{sgn}(\sigma)
    		\prod_{r=1}^{m}
    		a_{\sigma(2r-1),\sigma(2r)}.
    	\end{equation}
    	The Pfaffian of a matrix of odd order is defined to be \(0\), and
    	the Pfaffian of the empty matrix is defined to be \(1\). For an
    	even subset \(I\), we write \(A_I\) for the principal submatrix
    	with row and column set \(I\).

    	Let \(0\leq j\leq \ell\). For
    	\(0\leq\nu\leq\lfloor(j+1)/2\rfloor\), define the Pfaffian
    	gradings of the graph matrix \(S^j\) by
    	\begin{equation}\label{eq:psi-definition}
    		\Psi_{j,\nu}
    		=
    		\sum_{\substack{I\subseteq\{0,\ldots,j\}\\
    		|I|=2\nu}}
    		(-1)^{\operatorname{inv}(I,I^c)}
    		\operatorname{Pf}(S_I^j)\,f_{I^c}.
    	\end{equation}
    	Here \(I\) and its complement \(I^c\) are written in increasing
    	order, \(S_I^j\) is the principal submatrix with row and column set
    	\(I\), and \(f_{I^c}\) is the ordered wedge product of the vectors
    	\(f_i\) with \(i\in I^c\). The sign
    	\((-1)^{\operatorname{inv}(I,I^c)}\) is the sign of the shuffle
    	which moves the indices in \(I\) in front of those in \(I^c\).
    	We set \(\Psi_{j,\nu}=0\) for every integer \(\nu\) outside the
    	stated range.
    \end{definition}

    It follows directly from the definition that
    \[
    	\Psi_{j,\nu}\in
    	\Lambda^{j+1-2\nu}\mathbb{C}^n,
    	\qquad
    	\Psi_{j,0}=f_0\wedge\cdots\wedge f_j=:F_j.
    \]
    The gradings are initially local holomorphic representatives on
    the regular set. They have a common transition weight, proved in
    Lemma~\ref{lem:psi-transition}. Common zeros of \(F_j\) and poles
    of the other gradings are handled by rationality and independent
    saturation in \S\ref{subsec:pfaffian-polynomiality}.

    \begin{example}\label{ex:pfaffian-l1}
    	For \(\ell=1\), substituting the matrix \(S^1\) computed in
    	Example~\ref{ex:graph-l1} gives
    	\begin{align*}
    		\Psi_{0,0}&=f_0,\\
    		\Psi_{1,0}&=f_0\wedge f_1,\\
    		\Psi_{1,1}&=x_0\norm{f_0}^2.
    	\end{align*}
    \end{example}

    \begin{example}\label{ex:pfaffian-l2}
    	For \(\ell=2\), substituting the matrix \(S^2\) computed in
    	Example~\ref{ex:graph-l2} gives
    	\begin{align*}
    		\Psi_{0,0}&=f_0,\\
    		\Psi_{1,0}&=f_0\wedge f_1,
    		\qquad
    		\Psi_{1,1}=x_0\norm{f_0}^2,\\
    		\Psi_{2,0}&=f_0\wedge f_1\wedge f_2,\\
    		\Psi_{2,1}
    		&=x_0\norm{f_0}^2f_2\notag\\
    		&\qquad
    		-\left[
    		\partial\bigl(x_0\norm{f_0}^2\bigr)
    		-x_0\norm{f_0}^2
    		\partial\log\left(
    		\norm{f_0}^2\norm{f_1}^2\right)
    		\right]f_1\notag\\
    		&\qquad
    		+x_1\left(
    		\norm{f_1}^2+|x_0|^2\norm{f_0}^2\right)f_0.
    	\end{align*}
    	Following \citep{li2023conformal}, the grading \(\Psi_{2,1}\)
    	is denoted by \(g\) in that reference.
    \end{example}

    \subsubsection{Recurrence in the replacement index}

    For \(0\leq i\leq j\), let
    \(\iota_{t^{j+1}}\) be the interior product on
    \(\operatorname{span}\{f_0,\ldots,f_j\}\) determined by
    \(\iota_{t^{j+1}}(f_i)=t_i^{j+1}\), where \(t^{j+1}\) is the new
    column in Definition~\ref{def:S-and-w}. Thus
    \begin{equation}
    	\iota_{t^{j+1}}
    	(f_{i_1}\wedge\cdots\wedge f_{i_d})
    	=
    	\sum_{r=1}^{d}(-1)^{r-1}t_{i_r}^{j+1}
    	f_{i_1}\wedge\cdots\wedge
    	\widehat{f_{i_r}}\wedge\cdots\wedge f_{i_d}.
    \end{equation}
    This functional is defined on the old Frenet wedge products; its
    extension to the ambient space will not affect the formulas below.

    \begin{lemma}[Pfaffian recurrence]\label{lem:psi-recurrence}
    	For \(0\leq j<\ell\) and every admissible \(\nu\),
    	\begin{equation}\label{eq:psi-recurrence}
    		\boxed{\displaystyle
    		\Psi_{j+1,\nu}
    		=
    		\Psi_{j,\nu}\wedge f_{j+1}
    		+(-1)^j
    		\iota_{t^{j+1}}\Psi_{j,\nu-1}.}
    	\end{equation}
    \end{lemma}

    \begin{proof}
    	In the extension from \(S^j\) to \(S^{j+1}\), an even subset of
    	indices either omits the new index \(j+1\) or contains it. If it
    	omits \(j+1\), the Pfaffian is unchanged and the complementary
    	wedge product gains the factor \(f_{j+1}\). If it contains
    	\(j+1\), the new index is paired with some old index \(i\),
    	producing the coefficient \(t_i^{j+1}\) and lowering the
    	Pfaffian degree by one. Moving the new index past the remaining
    	\(j\) indices contributes the sign \((-1)^j\). These two cases
    	are precisely the two terms in \eqref{eq:psi-recurrence}.
    \end{proof}

    \subsubsection{Holomorphicity of the Pfaffian gradings}

    \begin{theorem}[Holomorphicity]\label{thm:psi-holomorphic}
    	On each graph-regular coordinate chart, for every
    	\(0\leq j\leq \ell\) and every admissible \(\nu\),
    	\begin{equation}\label{eq:psi-holomorphic}
    		\boxed{\bar\partial\Psi_{j,\nu}=0.}
    	\end{equation}
    \end{theorem}

    \begin{example}\label{ex:holomorphic-j1}
    	For \(j=1\), set \(\eta=x_0\norm{f_0}^2\). By
    	Example~\ref{ex:pfaffian-l1}, \(\Psi_{1,1}=\eta\). The connection
    	equation \eqref{eq:S-connection} gives \(\bar\partial\eta=0\), and
    	hence
    	\[
    		\bar\partial\Psi_{1,1}=0.
    	\]
    \end{example}

    \begin{example}\label{ex:holomorphic-j2}
    	For \(j=2\), write
    	\[
    		\eta=x_0\norm{f_0}^2,
    		\qquad
    		A=\partial\eta-\eta\partial\log\left(
    		\norm{f_0}^2\norm{f_1}^2\right),
    	\]
    	\[
    		B=x_1\left(
    		\norm{f_1}^2+|x_0|^2\norm{f_0}^2\right).
    	\]
    	By Example~\ref{ex:pfaffian-l2},
    	\begin{equation*}
    		\Psi_{2,1}=\eta f_2-Af_1+Bf_0.
    	\end{equation*}
    	The matrix connection gives
    	\[
    		\bar\partial\eta=0,
    		\qquad
    		\bar\partial A
    		=-\lambda_1\eta,
    		\qquad
    		\bar\partial B
    		=-\lambda_0A.
    	\]
    	Using the Frenet equations, we obtain
    	\begin{align*}
    		\bar\partial(\eta f_2)
    		&=-\lambda_1\eta f_1,\\
    		\bar\partial(-Af_1)
    		&=\lambda_1\eta f_1+\lambda_0Af_0,\\
    		\bar\partial(Bf_0)
    		&=-\lambda_0Af_0.
    	\end{align*}
    	Therefore \(\bar\partial\Psi_{2,1}=0\).
    \end{example}

    \begin{proof}
        Fix \(j,\nu\) and work on a graph-regular coordinate chart.
        The Frenet vectors \(f_0,\ldots,f_j\) are linearly independent;
        write \(W_j=\operatorname{span}\{f_0,\ldots,f_j\}\).
        All contractions below act directly on \(\Lambda W_j\), so no
        auxiliary vector space or comparison map is needed.

        \smallskip
        \noindent\emph{1. Contractions in the Frenet span.}
        For \(0\leq i\leq j\), let \(\iota_i\) be contraction by the
        linear functional on \(W_j\) which takes \(f_k\) to
        \(\delta_{ik}\). Thus, on an ordered wedge of Frenet vectors,
        \(\iota_i\) removes \(f_i\), with the sign determined by its
        position, and gives zero if \(f_i\) does not occur. Set
        \begin{equation}
            \mathcal I_{S^j}
            =\sum_{0\leq i<k\leq j}
              S_{ik}^j\,\iota_k\iota_i,
            \qquad
            \operatorname{vol}_j=f_0\wedge\cdots\wedge f_j.
        \end{equation}
        For an even set \(I=\{i_1<\cdots<i_{2\nu}\}\), the contraction
        removing precisely these indices satisfies
        \[
            \iota_{i_{2\nu}}\cdots\iota_{i_1}
            \operatorname{vol}_j
            =(-1)^{\operatorname{inv}(I,I^c)}f_{I^c}.
        \]
        In \(\mathcal I_{S^j}^{\nu}\), the selected \(2\nu\)
        indices are paired in every possible way. The signs of the
        pairings and their products of matrix entries give
        \(\operatorname{Pf}(S_I^j)\), by
        \eqref{eq:matrix-pfaffian}; the \(\nu\) pair-contraction
        operators commute with one another, so each unordered
        collection of pairs occurs \(\nu!\) times. Consequently,
        the defining sum \eqref{eq:psi-definition} is exactly
        \begin{equation}\label{eq:psi-contraction}
            \Psi_{j,\nu}
            =\frac{1}{\nu!}\,
              \mathcal I_{S^j}^{\nu}\operatorname{vol}_j.
        \end{equation}
        This includes \(\nu=0\), for which the right-hand side
        is \(\operatorname{vol}_j\).

        \smallskip
        \noindent\emph{2. Equivariance.}
        Let \(g\in GL(j+1,\mathbb C)\), and put
        \(f'_i=\sum_u g_{iu}f_u\) and
        \(S'=gS^jg^{\transp}\). The new tuple spans the same
        space \(W_j\). Denote by \(\iota'_i\) the contractions
        characterized by \(\iota'_i(f'_k)=\delta_{ik}\).
        Expressing each old vector in the new tuple gives
        \[
            \iota'_i
            =\sum_u(g^{-1})_{ui}\iota_u.
        \]
        Since contractions anticommute, we may write
        \[
            \mathcal I_{S^j}
            =\frac12\sum_{i,k=0}^j
              S_{ik}^j\,\iota_k\iota_i.
        \]
        Substitution of \(S'=gS^jg^{\transp}\) and the formula
        for \(\iota'_i\) therefore gives, as operators on
        \(\Lambda W_j\),
        \[
            \mathcal I_{S'}
            =\frac12\sum_{u,v}
             \bigl(g^{-1}S'(g^{-1})^{\transp}\bigr)_{uv}
             \iota_v\iota_u
            =\mathcal I_{S^j}.
        \]
        On the other hand,
        \(f'_0\wedge\cdots\wedge f'_j
        =\det(g)\operatorname{vol}_j\).
        Applying \eqref{eq:psi-contraction} to the primed data,
        we obtain
        \begin{equation}\label{eq:psi-equivariance}
            \Psi_{j,\nu}(gf,gS^jg^{\transp})
            =\det(g)\Psi_{j,\nu}(f,S^j).
        \end{equation}
        Here \(gf\) means the ordered tuple
        \(f'_i=\sum_u g_{iu}f_u\); the identity holds for every
        invertible \(g\), independently of any differential equation.

        \smallskip
        \noindent\emph{3. Differentiation along the connection.}
        Let \(F=(f_0,\ldots,f_j)^{\transp}\), and define
        the strictly lower triangular matrix \(\mathsf L\) by
        \begin{equation}
            \mathsf L_{i,i-1}=\lambda_{i-1},
            \qquad 1\leq i\leq j,
        \end{equation}
        with every other entry zero. The Frenet equations and
        the already established connection equation
        \eqref{eq:S-connection} say, respectively,
        \begin{equation}
            \bar\partial F=-\mathsf L F,
        \end{equation}
        \begin{equation}
            \bar\partial S^j
            =-\mathsf L S^j-S^j\mathsf L^{\transp}.
        \end{equation}

        To justify the differentiation, regard the right-hand
        side of \eqref{eq:psi-definition} as a function
        \(\Phi(F,S)\) of an ordered vector tuple and a
        skew-symmetric matrix. It is polynomial in their complex
        entries and contains no conjugates. Thus its complex
        differential \(D\Phi\) obeys the \(\bar\partial\) chain
        rule, giving at each fixed point of the chart
        \[
            \bar\partial\Psi_{j,\nu}
            =D\Phi_{(F,S^j)}
             \bigl(-\mathsf L F,\,
                   -\mathsf L S^j-S^j\mathsf L^{\transp}\bigr).
        \]
        Freeze \(F,S^j,\mathsf L\) at this point and introduce
        the \emph{auxiliary} complex parameter \(\varepsilon\) and
        \(g_\varepsilon=I-\varepsilon\mathsf L\). Direct expansion gives
        \[
            \begin{aligned}
              g_\varepsilon F
                &=F-\varepsilon\mathsf L F,\\
              g_\varepsilon S^jg_\varepsilon^{\transp}
                &=S^j-\varepsilon
                  (\mathsf L S^j+S^j\mathsf L^{\transp})
                  +\varepsilon^2\mathsf L S^j\mathsf L^{\transp}.
            \end{aligned}
        \]
        The two inputs of \(\Phi\) therefore have precisely the
        first-order increments appearing in the chain rule. By
        the definition of the differential,
        \[
            \bar\partial\Psi_{j,\nu}
            =\left.\frac{d}{d\varepsilon}\right|_{\varepsilon=0}
              \Phi(g_\varepsilon F,\,
                   g_\varepsilon S^jg_\varepsilon^{\transp}).
        \]
        Finally, \eqref{eq:psi-equivariance} evaluates the last
        derivative. Because \(g_\varepsilon\) is unit lower
        triangular, \(\det(g_\varepsilon)=1\) for every
        \(\varepsilon\), and \(\operatorname{tr}\mathsf L=0\).
        Hence
        \begin{equation}
            \bar\partial\Psi_{j,\nu}
            =\left.\frac{d}{d\varepsilon}\right|_{\varepsilon=0}
               \det(I-\varepsilon\mathsf L)\Psi_{j,\nu}
            =-(\operatorname{tr}\mathsf L)\Psi_{j,\nu}=0.
        \end{equation}
        This proves holomorphicity using only the Frenet equations
        and the graph connection from
        Lemma~\ref{lem:simultaneous}.
    \end{proof}

    \subsubsection{Decomposition of the norm}

    \begin{lemma}[Frenet form of the Pfaffian norm]
    	\label{lem:pfaffian-norm-expansion}
    	For every \(0\leq j\leq \ell\),
    	\begin{equation}\label{eq:psi-norm-expanded}
    		\sum_{\nu=0}^{\lfloor(j+1)/2\rfloor}
    		\norm{\Psi_{j,\nu}}^2
    		=
    		\sum_{\substack{I\subseteq\{0,\ldots,j\}\\
    		|I|\ \text{even}}}
    		|\operatorname{Pf}(S_I^j)|^2
    		\prod_{k\notin I}\norm{f_k}^2.
    	\end{equation}
    \end{lemma}

    \begin{proof}
    	The vectors \(f_i\) are Hermitian orthogonal. Therefore the wedge
    	products associated with distinct Frenet index sets are mutually
    	orthogonal, and \eqref{eq:psi-definition} can be expanded in
    	squared norms without cross terms. This gives
    	\eqref{eq:psi-norm-expanded}.
    \end{proof}

    \begin{lemma}[Pure Pfaffian identity]
    	\label{lem:pure-pfaffian-identity}
    	Let \(K\) be a \(d\times d\) complex skew-symmetric matrix. Then
    	\begin{equation}\label{eq:pure-skew-identity}
    		\left(
    		\sum_{I\ \text{even}}
    		|\operatorname{Pf}(K_I)|^2
    		\right)^2
    		=
    		\det(I+KK^*).
    	\end{equation}
    \end{lemma}

    \begin{proof}
    	Let \(K^*\) denote the conjugate transpose. Consider the
    	\(2d\times2d\) skew-symmetric block matrix
    	\begin{equation}
    		\mathcal B(K)=
    		\begin{pmatrix}
    			K&I\\
    			-I&-\overline K
    		\end{pmatrix}.
    	\end{equation}
    	In the expansion of its Pfaffian, taking an even set \(I\) of
    	mutually paired indices in the upper block forces the same index
    	set to be paired in the lower block. All remaining indices are
    	joined by the identity blocks. The two halves contribute
    	\(\operatorname{Pf}(K_I)\) and
    	\(\overline{\operatorname{Pf}(K_I)}\), and the common shuffle
    	sign is \((-1)^{d(d-1)/2}\). Hence
    	\begin{equation}
    		\operatorname{Pf}\mathcal B(K)
    		=
    		(-1)^{d(d-1)/2}
    		\sum_{I\ \text{even}}
    		|\operatorname{Pf}(K_I)|^2.
    	\end{equation}
    	After interchanging block columns and applying Schur elimination,
    	one obtains
    	\begin{equation}
    		\det\mathcal B(K)
    		=\det(I-K\overline K)
    		=\det(I+KK^*),
    	\end{equation}
    	where \(K^*=-\overline K\) has been used. Squaring the preceding
    	two formulas and using
    	\((\operatorname{Pf}B)^2=\det B\) proves
    	\eqref{eq:pure-skew-identity}.
    \end{proof}

    \begin{lemma}[Single-step growth]\label{lem:pfaffian-growth}
    	Let
    	\begin{equation}\label{eq:normalized-K}
    		K_{ik}
    		=\frac{S_{ik}^j}
    		{\norm{f_i}\norm{f_k}},
    	\end{equation}
    	and let
    	\begin{equation}
    		K'=
    		\begin{pmatrix}
    			K&t\\
    			-t^{\transp}&0
    		\end{pmatrix},
    		\qquad
    		t_i
    		=\frac{t_i^{j+1}}
    		{\norm{f_i}\norm{f_{j+1}}}.
    	\end{equation}
    	Then
    	\begin{equation}\label{eq:pf-growth-ratio}
    		\frac{
    		\sum_{I\ \text{even}}
    		|\operatorname{Pf}(K'_I)|^2}
    		{
    		\sum_{I\ \text{even}}
    		|\operatorname{Pf}(K_I)|^2}
    		=
    		1+t^*(I+KK^*)^{-1}t.
    	\end{equation}
    \end{lemma}

    \begin{proof}
    	Write \(M=I+KK^*\), which is positive definite Hermitian. The
    	new indices in \(\mathcal B(K')\) may be moved to the end,
    	producing a matrix of the form
    	\begin{equation}
    		\begin{pmatrix}
    			\mathcal B(K)&C\\
    			-C^{\transp}&J
    		\end{pmatrix},
    		\qquad
    		C=
    		\begin{pmatrix}
    			t&0\\
    			0&-\overline t
    		\end{pmatrix},
    		\quad
    		J=
    		\begin{pmatrix}
    			0&1\\
    			-1&0
    		\end{pmatrix}.
    	\end{equation}
    	Solving the two linear systems defining
    	\(\mathcal B(K)^{-1}\) gives
    	\begin{equation}
    		\mathcal B(K)^{-1}
    		=
    		\begin{pmatrix}
    			-\overline K M^{-1}&-(I+K^*K)^{-1}\\
    			M^{-1}&M^{-1}K
    		\end{pmatrix}.
    	\end{equation}
    	Substitution into the Schur complement
    	\(J+C^{\transp}\mathcal B(K)^{-1}C\) shows that its two
    	diagonal entries vanish by skew-symmetry and that its upper-right
    	entry is the real number \(1+t^*M^{-1}t\). Hence
    	\begin{equation}
    		J+C^{\transp}\mathcal B(K)^{-1}C
    		=(1+t^*M^{-1}t)J.
    	\end{equation}
    	Simultaneous row and column operations on a skew-symmetric block
    	matrix multiply the Pfaffian by the corresponding determinant
    	factors. After accounting for the reordering signs, this gives
    	\eqref{eq:pf-growth-ratio}.
    \end{proof}

    \begin{theorem}[Norm decomposition]\label{thm:pfaffian-norm}
    	For every \(0\leq j\leq \ell\),
    	\begin{equation}\label{eq:psi-norm-identity}
    		\boxed{\displaystyle
    		\sum_{\nu=0}^{\lfloor(j+1)/2\rfloor}
    		\norm{\Psi_{j,\nu}}^2
    		=
    		\prod_{i=0}^{j}\norm{v_i}^2.}
    	\end{equation}
    	In particular,
    	\begin{equation}\label{eq:final-psi-norm-identity}
    		\sum_{\nu=0}^{\lfloor(\ell+1)/2\rfloor}
    		\norm{\Psi_{\ell,\nu}}^2
    		=
    		\prod_{i=0}^{\ell}\norm{v_i}^2.
    	\end{equation}
    \end{theorem}

    \begin{example}\label{ex:norm-decomposition-j0}
    	For \(j=0\), we have \(\Psi_{0,0}=f_0\) and \(v_0=f_0\).
    	Consequently,
    	\[
    		\sum_{\nu=0}^{0}\norm{\Psi_{0,\nu}}^2
    		=\norm{\Psi_{0,0}}^2
    		=\norm{f_0}^2
    		=\norm{v_0}^2.
    	\]
    \end{example}

    \begin{example}\label{ex:norm-decomposition-j1}
    	For \(j=1\), we have \(v_0=f_0\),
    	\(v_1=f_1+x_0\overline{f_0}\),
    	\(\Psi_{1,0}=f_0\wedge f_1\), and
    	\(\Psi_{1,1}=x_0\norm{f_0}^2\). Hence
    	\begin{align*}
    		\norm{\Psi_{1,0}}^2+|\Psi_{1,1}|^2
    		&=\norm{f_0}^2\norm{f_1}^2
    		+|x_0|^2\norm{f_0}^4\\
    		&=\norm{v_0}^2
    		\left(\norm{f_1}^2
    		+|x_0|^2\norm{f_0}^2\right)
    		=\norm{v_0}^2\norm{v_1}^2.
    	\end{align*}
    \end{example}

    \begin{proof}
    	The cases \(j=0\) and \(j=1\) are the two preceding examples.
    	We prove the induction step from
    	level \(j\) to level \(j+1\). Normalize the graph matrix at
    	level \(j\) as in \eqref{eq:normalized-K}. Multi-linearity
    	of the Pfaffian gives, for each even subset \(I\),
    	\begin{equation}
    		\operatorname{Pf}(S_I^j)
    		=
    		\left(\prod_{i\in I}\norm{f_i}\right)
    		\operatorname{Pf}(K_I).
    	\end{equation}
    	Consequently,
    	\eqref{eq:psi-norm-expanded} becomes
    	\begin{equation}\label{eq:psi-norm-normalized}
    		\sum_{\nu}
    		\norm{\Psi_{j,\nu}}^2
    		=
    		\left(\prod_{i=0}^{j}\norm{f_i}^2\right)
    		\sum_{I\ \text{even}}
    		|\operatorname{Pf}(K_I)|^2.
    	\end{equation}
    	By Lemma~\ref{lem:pfaffian-growth}, passing to \(j+1\) multiplies
    	the sum in the last line by
    	\(1+t^*(I+KK^*)^{-1}t\).

    	It remains to identify this factor with the norm of \(v_{j+1}\).
    	Set \(e_i=f_i/\norm{f_i}\) and write
    	\begin{equation}
    		\frac{v_{j+1}}{\norm{f_{j+1}}}
    		=e_{j+1}
    		+\sum_{i=0}^{j}
    		\left(
    		\alpha_i e_i+\gamma_i\overline{e_i}
    		\right),
    	\end{equation}
    	where
    	\begin{equation}
    		\alpha_i
    		=
    		\frac{p_i^{j+1}\norm{f_i}}
    		{\norm{f_{j+1}}},
    		\qquad
    		\gamma_i
    		=
    		\frac{q_i^{j+1}\norm{f_i}}
    		{\norm{f_{j+1}}}.
    	\end{equation}
    	The normalized old graph basis is
    	\(e_i+\sum_kK_{ki}\overline{e_k}\). Since the new vector is
    	Hermitian orthogonal to the old graph, imposing the single
    	vector equation
    	\(\alpha+K^*\gamma=0\) is equivalent to testing
    	orthogonality against every graph basis vector. Thus
    	\begin{equation}
    		\alpha=-K^*\gamma.
    	\end{equation}
    	The definition of the normalized new column has the
    	\(e_i\)-components
    	\(t_i=\gamma_i-\sum_kK_{ik}\alpha_k\), or in vector form
    	\(t=\gamma-K\alpha\). Substitution therefore yields
    	\begin{equation}
    		t=\gamma-K\alpha=(I+KK^*)\gamma.
    	\end{equation}
    	As \(I+KK^*\) is positive definite, this identity determines
    	\(\gamma\) and gives
    	\begin{align}
    		\frac{\norm{v_{j+1}}^2}{\norm{f_{j+1}}^2}
    		&=1+\alpha^*\alpha+\gamma^*\gamma\notag\\
    		&=1+\gamma^*(I+KK^*)\gamma
    		=1+t^*(I+KK^*)^{-1}t.
    		\label{eq:replacement-norm-ratio}
    	\end{align}
    	Thus the growth factor for the normalized Pfaffian sum is
    	\(\norm{v_{j+1}}^2/\norm{f_{j+1}}^2\).

    	The prefactor in \eqref{eq:psi-norm-normalized} acquires an
    	additional factor \(\norm{f_{j+1}}^2\) when \(j\) is replaced by
    	\(j+1\). Multiplying this with the growth factor gives
    	\(\norm{v_{j+1}}^2\). The induction therefore proves
    	\eqref{eq:psi-norm-identity} at every level. Taking \(j=\ell\)
    	yields \eqref{eq:final-psi-norm-identity}.
    \end{proof}

\subsection{Covariance and the common ambient space}

We record the covariance on the common regular set.

    \begin{lemma}[Common transition weight]
    	\label{lem:psi-transition}
    	Under a holomorphic change of affine coordinate
    	\(z=z(\zeta)\), write \(\tau=dz/d\zeta\). If
    	\(f_0^{\mathrm{new}}=\chi f_0\) with \(\chi\neq0\), then
    	\begin{align}
    		f_i^{\mathrm{new}}
    		&=\chi\tau^if_i,
    		&
    		v_j^{\mathrm{new}}
    		&=\chi\tau^jv_j,\label{eq:transition-frenet-v}\\
    		(S_{ik}^j)^{\mathrm{new}}
    		&=\chi^2\tau^{i+k}S_{ik}^j,
    		&
    		\Psi_{j,\nu}^{\mathrm{new}}
    		&=\chi^{j+1}
    		\tau^{j(j+1)/2}\Psi_{j,\nu}.\label{eq:transition-psi}
    	\end{align}
    	For \(0\leq j<\ell\), the deviation parameter transforms by
    	\begin{equation}\label{eq:transition-deviation}
    		(x_j)^{\mathrm{new}}
    		=\frac{\chi}{\overline\chi}\,
    		\tau^{j+1}\overline\tau^{-j}x_j.
    	\end{equation}
    \end{lemma}

    \begin{proof}
    	Frenet orthogonalization removes the tangential component at each
    	step, so the first identity follows by induction on \(i\). The
    	retained line is unchanged geometrically and its highest
    	\(f_j\)-coefficient remains one, which gives the identity for
    	\(v_j\).

    	To derive \eqref{eq:transition-deviation}, note that
    	\((u_j)^{\mathrm{new}}=\chi\tau^{j+1}u_j\). Combining this
    	with
    	\((\overline{v_j})^{\mathrm{new}}
    	=\overline\chi\,\overline\tau^{\,j}\overline{v_j}\)
    	in \(v_{j+1}=u_j+x_j\overline{v_j}\) and comparing the
    	\(\overline{v_j}\)-components gives
    	\eqref{eq:transition-deviation}. The multiplier is nonzero on
    	the common regular chart, so the vanishing of \(x_j\) is
    	independent of the chosen coordinate and lift.

    	The graph matrix \(S^j\) has weight
    	\(\chi^2\tau^{i+k}\) in the \((i,k)\)-entry. In a summand of
    	\(\Psi_{j,\nu}\), every selected index from the Pfaffian
    	contributes \(\chi^2\tau^{i+k}\), while every index in the
    	complement contributes \(\chi\tau^k\). Since \(I\) and \(I^c\)
    	partition \(\{0,\ldots,j\}\), the total factor is
    	\begin{equation}
    		\chi^{j+1}
    		\tau^{0+1+\cdots+j}
    		=
    		\chi^{j+1}\tau^{j(j+1)/2},
    	\end{equation}
    	independently of \(\nu\). This proves
    	\eqref{eq:transition-psi}.
    \end{proof}

Place all gradings in the fixed Hermitian space
\[
 \mathcal V=\bigoplus_{k=0}^{\ell+1}\Lambda^k\mathbb C^n,
 \qquad
 \boldsymbol\Psi_j=\bigoplus_\nu\Psi_{j,\nu},\quad
 \boldsymbol\Psi_{-1}=1,\quad P_{-1}=1.
\]
Missing exterior degrees are filled by zero. Thus
Theorem~\ref{thm:pfaffian-norm} reads
\begin{equation}\label{eq:raw-augmented-norm}
 \norm{\boldsymbol\Psi_j}^2=P_j=\prod_{i=0}^j\norm{v_i}^2.
\end{equation}
The common ambient space is also necessary for the mixed kernels:
one must pair every matching exterior degree, not just the top wedges.

\subsection{Actual energies and exact raw equations}

Define
\begin{equation}\label{eq:rigidity-sigma}
 \sigma_j=\frac{\norm{u_j}^2}{\norm{v_j}^2},\qquad
 0\leq j\leq\ell,\qquad \sigma_{-1}=0.
\end{equation}
Here \(u_\ell\) is defined by \eqref{eq:deviation-forward} without
performing an additional replacement. The next identities depend only
on the short-window graph calculus.

    \begin{lemma}[Connection and metric identities]
    	\label{lem:rigidity-connection}
    	Under the range in Corollary~\ref{cor:bew-admissible-range},
    	\begin{align}
    		\bar\partial u_j
    		&=-\sigma_jv_j+c_j\overline{v_j},
    		&&0\leq j\leq \ell,\label{eq:rigidity-conn1}\\
    		\bar\partial v_{j+1}
    		&=-\sigma_jv_j+x_j\overline{u_j},
    		&&0\leq j<\ell,\label{eq:rigidity-conn2}\\
    		\partial\bar\partial\log\norm{v_j}^2
    		&=\sigma_j-\sigma_{j-1},
    		&
    		\partial\bar\partial\log P_j&=\sigma_j,
    		\label{eq:rigidity-curv-pot}\\
    		\partial x_{j-1}
    		-\bigl(\partial\log\norm{v_j}^2\bigr)x_{j-1}
    		&=c_j,
    		&&1\leq j<\ell,\label{eq:rigidity-src1}
    	\end{align}
    	Moreover, for \(0\leq j<\ell\),
    	\begin{equation}
    		\bar\partial x_j
    		+\bigl(\bar\partial\log\norm{v_j}^2\bigr)x_j
    		=-c_j,
    	\end{equation}
    	and, for \(0\leq j\leq\ell\),
    	\begin{equation}\label{eq:rigidity-metric}
    		g_j
    		=
    		2(\sigma_j+\sigma_{j-1})
    		\,\mathrm{d}z\,\mathrm{d}\overline{z},
    	\end{equation}
    	\begin{equation}
    		K_j
    		=
    		-\frac{\partial\bar\partial\log(\sigma_j+\sigma_{j-1})}
    		{\sigma_j+\sigma_{j-1}}.
    	\end{equation}
    \end{lemma}

    \begin{proof}
    	The decompositions \eqref{eq:Ubar-decomposition} and
    	\eqref{eq:terminal-graph-connection} give
    	\(\bar\partial u_j=Z_j+c_j\overline{v_j}\) with
    	\(Z_j\in\mathbb L_j\). Differentiate
    	\(H(u_j,w_i^j)=0\) for \(i<j\). The contact condition
    	\eqref{eq:contact-condition} and
    	\(b(u_j,v_j)=0\) show that \(Z_j\) is Hermitian orthogonal to
    	\(w_i^j\) for every \(i<j\). Hence \(Z_j=\gamma_jv_j\).
    	Differentiating \(H(u_j,v_j)=0\), using
    	\[
    		\partial v_j
    		=
    		u_j+
    		\bigl(\partial\log\norm{v_j}^2\bigr)v_j,
    	\]
    	gives
    	\(\gamma_j\norm{v_j}^2+\norm{u_j}^2=0\). This proves
    	\eqref{eq:rigidity-conn1}.

    	Apply \(\bar\partial\) to
    	\(v_{j+1}=u_j+x_j\overline{v_j}\). The
    	\(\overline{v_j}\)-terms cancel by
    	\eqref{eq:scalar-replacement}, yielding
    	\eqref{eq:rigidity-conn2}. Since
    	\(u_j\perp_H\overline{v_j}\),
    	\begin{equation}\label{eq:rigidity-pivot}
    		\frac{\norm{v_{j+1}}^2}{\norm{v_j}^2}
    		=
    		\sigma_j
    		+|x_j|^2.
    	\end{equation}
    	The two terms in \eqref{eq:rigidity-conn2} are orthogonal.
    	Consequently,
    	\[
    		\frac{\norm{\bar\partial v_{j+1}}^2}
    		{\norm{v_{j+1}}^2}
    		=
    		\sigma_j.
    	\]
    	For completeness, put \(h_j=\norm{v_j}^2\). The graph
    	connection and \eqref{eq:orthogonality-v-w} give
    	\[
    		H(\bar\partial v_j,v_j)=0.
    	\]
    	With \(H\) linear in its first argument,
    	\[
    		\partial h_j=H(\partial v_j,v_j),
    		\qquad
    		H(\partial\bar\partial v_j,v_j)
    		=-\norm{\bar\partial v_j}^2,
    	\]
    	where the second identity follows by differentiating
    	\(H(\bar\partial v_j,v_j)=0\). Hence
    	\[
    		\partial\bar\partial h_j
    		=\norm{\partial v_j}^2-\norm{\bar\partial v_j}^2.
    	\]
    	By \eqref{eq:deviation-forward} and orthogonal projection,
    	\[
    		\norm{u_j}^2
    		=\norm{\partial v_j}^2
    		-\frac{|H(\partial v_j,v_j)|^2}{h_j}.
    	\]
    	Therefore
    	\[
    		\partial\bar\partial\log h_j
    		=\frac{\norm{u_j}^2}{h_j}
    		-\frac{\norm{\bar\partial v_j}^2}{h_j}
    		=\sigma_j-\sigma_{j-1}.
    	\]
    	The case \(j=0\) uses \(\bar\partial v_0=0\) and
    	\(\sigma_{-1}=0\). This proves the second identity in
    	\eqref{eq:rigidity-curv-pot}.

    	For \(j\geq1\), compare the coefficient of
    	\(\overline{f_j}\) in
    	\[
    		\bar\partial u_j
    		=
    		\bar\partial\left(
    		\partial v_j-
    		\bigl(\partial\log\norm{v_j}^2\bigr)v_j
    		\right).
    	\]
    	Indeed, \eqref{eq:q-recurrence} gives
    	\(q_{j-1}^j=x_{j-1}\), so the coefficient of
    	\(\overline{f_{j-1}}\) in \(u_j\) is
    	\[
    	 \partial x_{j-1}-(\partial\log\norm{v_j}^2)x_{j-1}.
    	\]
    	The only contribution to \(\overline{f_j}\) after
    	\(\bar\partial\)-differentiation comes from this coefficient.
    	On the right of \eqref{eq:rigidity-conn1}, the coefficient is
    	\(c_j\), proving \eqref{eq:rigidity-src1}. The remaining scalar
    	equation is \eqref{eq:scalar-replacement}, for \(0\leq j<\ell\).
    	Finally, the two
    	normal energies of the real layer
    	\(V_j=v_j\oplus\overline{v_j}\) are \(\sigma_j\) and
    	\(\sigma_{j-1}\). The stated metric and curvature formulas follow
    	from the normalization of the Grassmannian metric and the usual
    	formula for a conformal metric.
    \end{proof}

In particular, \(\sigma_j\,dz\,d\bar z\) is the true forward-energy
form of the smooth line map \(\underline v_j\), and
\(\sigma_{j-1}\,dz\,d\bar z\) is its backward-energy form. This can
also be read directly from the two normal components of
\(dv_j\). On changing to any smooth nonzero frame of that actual
line, these energies are unchanged. They therefore extend smoothly
and nonnegatively across every excluded point.

Define, still in the raw gauge,
\begin{equation}\label{eq:raw-mixed-definition}
 B_j=H(\boldsymbol\Psi_{j+1},\boldsymbol\Psi_{j-1}),\qquad0\leq j<\ell.
\end{equation}

    \begin{lemma}[Pfaffian mixed-kernel identity]
    	\label{lem:rigidity-mixed}
    	For \(0\leq j<\ell\),
    	\begin{equation}\label{eq:rigidity-mixed}
    		\boxed{B_j=x_jP_j}.
    	\end{equation}
    \end{lemma}

    \begin{proof}
    	Let \(K\) be a skew-symmetric matrix, \(M=I+KK^*\), and
    	\begin{equation}
    		\mathscr N(K)=
    		\sum_{I\ \text{even}}
    		|\operatorname{Pf}(K_I)|^2.
    	\end{equation}
    	For vectors \(\mathsf v,\mathsf t\) and a scalar \(\omega\),
    	the bordered-Pfaffian identity is
    	\begin{equation}\label{eq:rigidity-pf-border}
    		\begin{aligned}
    		&\sum_{I\ \text{even}}
    		\operatorname{Pf}
    		\begin{pmatrix}
    			K_I&\mathsf v_I&\mathsf t_I\\
    			-\mathsf v_I^{\transp}&0&\omega\\
    			-\mathsf t_I^{\transp}&-\omega&0
    		\end{pmatrix}
    		\overline{\operatorname{Pf}(K_I)}\\
    		&\qquad
    		=
    		\mathscr N(K)
    		\left(
    		\omega-\mathsf v^{\transp}\overline K M^{-1}\mathsf t
    		\right).
    		\end{aligned}
    	\end{equation}
    	Appendix~\ref{app:bordered-pfaffian} proves
    	\eqref{eq:rigidity-pf-border} by block expansion and a
    	Pfaffian Schur complement, including the permutation signs
    	and the required inverse block.

    	Now use the orthonormal frame
    	\(e_i=f_i/\norm{f_i}\). Let \(K\) be the correspondingly
    	normalized principal matrix for the old indices
    	\(0,\ldots,j-1\), and let \(\mathsf v,\mathsf t\) be the
    	normalized columns at indices \(j\) and \(j+1\). Put
    	\begin{equation}
    		\omega
    		=
    		\frac{S_{j,j+1}}
    		{\norm{f_j}\norm{f_{j+1}}},
    		\qquad
    		\frac{v_j}{\norm{f_j}}
    		=
    		e_j+\sum_{i<j}
    		(\alpha_i e_i+\gamma_i\overline{e_i}).
    	\end{equation}
    	Graph membership and orthogonality to the old graph give
    	\begin{equation}
    		\alpha=\overline K\gamma,
    		\qquad
    		\mathsf v=\gamma-K\alpha=M\gamma,
    		\qquad
    		\alpha=\overline K M^{-1}\mathsf v.
    	\end{equation}
    	Since \(\overline K M^{-1}\) is skew-symmetric,
    	\(\omega+\alpha^{\transp}\mathsf t
    	=\omega-\mathsf v^{\transp}\overline K M^{-1}\mathsf t\).

    	It remains to identify the exterior degrees and restore the
    	normalization factors. The summand
    	\(\Psi_{j+1,\nu'}\) has degree \(j+2-2\nu'\), while
    	\(\Psi_{j-1,\nu''}\) has degree \(j-2\nu''\). A nonzero pairing
    	between the two augmented vectors therefore requires
    	\(\nu'=\nu''+1\). For an even set \(I\) of size \(2\nu''\),
    	the normalized Pfaffians and the complementary wedge factors
    	contain, respectively,
    	\(\prod_{i\in I}\norm{f_i}^2\) and
    	\(\prod_{i<j,\,i\notin I}\norm{f_i}^2\). Their product is
    	\[
    		\prod_{i=0}^{j-1}\norm{f_i}^2=P_{j-1}.
    	\]
    	The remaining two columns are exactly
    	\(\eta_i=p_i^j\) and \(S_{i,j+1}\). Equivalently, using
    	\(P_{j-1}=(\prod_{i<j}\norm{f_i}^2)\mathscr N(K)\) and
    	\[
    		\eta_i
    		=\sqrt{\norm{f_j}^2/\norm{f_i}^2}\,\alpha_i
    		=p_i^j\quad(i<j),
    		\qquad \eta_j=p_j^j=1,
    	\]
    	\eqref{eq:rigidity-pf-border} gives
    	\begin{equation}
    		H(\boldsymbol\Psi_{j+1},\boldsymbol\Psi_{j-1})
    		=
    		P_{j-1}\sum_{i=0}^{j}\eta_iS_{i,j+1}
    		=P_{j-1}x_j\norm{v_j}^2
    		=x_jP_j.
    	\end{equation}
    	For the second equality, the new column satisfies
    	\(S_{i,j+1}=b(v_{j+1},w_i^j)\). Therefore
    	\[
    	 \sum_{i=0}^{j}\eta_iS_{i,j+1}
    	 =b(v_{j+1},v_j)=x_j\norm{v_j}^2,
    	\]
    	by \eqref{eq:v-in-w}, \(v_{j+1}=u_j+x_j\overline{v_j}\),
    	and \(b(u_j,v_j)=0\).
    	This is an identity in the original monic gauge on the regular set.
    \end{proof}

\begin{proposition}[Raw norm and source equations]
\label{prop:rigidity-bridge}
The raw data satisfy
\begin{align}
 P_{j-1}P_{j+1}
 &=P_j^2\partial\bar\partial\log P_j+B_jB_j^*,
 &&0\leq j<\ell,\label{eq:rigidity-N}\\
 P_j\partial B_{j-1}-(\partial P_j)B_{j-1}
 &=-P_{j-1}\bar\partial B_j+(\bar\partial P_{j-1})B_j,
 &&1\leq j<\ell.\label{eq:rigidity-C}
\end{align}
On the real diagonal \(B_j^*=\overline{B_j}\); after polarization
the star interchanges the two variables and conjugates coefficients.
In particular, \(B_0\) is the scalar grading of \(\boldsymbol\Psi_1\).
\end{proposition}

\begin{proof}
By \eqref{eq:rigidity-pivot},
\(P_{j-1}P_{j+1}/P_j^2=\sigma_j+|x_j|^2\).
Substitute \eqref{eq:rigidity-curv-pot} and
\eqref{eq:rigidity-mixed} to obtain \eqref{eq:rigidity-N}.
For the source equation, differentiating \(B_{j-1}=x_{j-1}P_{j-1}\)
and \(B_j=x_jP_j\) gives
\begin{align*}
 P_j\partial B_{j-1}-(\partial P_j)B_{j-1}
 &=P_jP_{j-1}\bigl(\partial x_{j-1}
       -(\partial\log\norm{v_j}^2)x_{j-1}\bigr),\\
 -P_{j-1}\bar\partial B_j+(\bar\partial P_{j-1})B_j
 &=-P_jP_{j-1}\bigl(\bar\partial x_j
       +(\bar\partial\log\norm{v_j}^2)x_j\bigr).
\end{align*}
Both expressions equal \(P_jP_{j-1}c_j\) by
Lemma~\ref{lem:rigidity-connection}. This proves the exact source
equation, including its sign. The assertion about \(B_0\) follows
from \(\boldsymbol\Psi_{-1}=1\).
\end{proof}

The rigidity argument uses both the norm and source equations in this
original gauge.

\subsection{Rationality and independent saturation}
\label{subsec:pfaffian-polynomiality}

\begin{lemma}[Rational holomorphic representatives]
\label{lem:rational-pfaffian-data}
For a primitive polynomial seed lift, every ambient component of
\(\boldsymbol\Psi_j\) belongs to \(\mathbb C(z)\). The identities
\eqref{eq:rigidity-N}--\eqref{eq:rigidity-C} therefore hold as
rational identities after polarization in \(\mathbb C(z,w)\).
\end{lemma}

\begin{proof}
The seed Frenet vectors are rational in \(z,\bar z\), by their
orthogonalization formulas. Start with the canonical rational
representative \(A_0\) of the seed. Proposition~5.1 of
\citep{bahy1989explicit} states that the components of every
globally defined meromorphic section of a harmonic bundle are
rational in \(z,\bar z\). Lemma~5.2 then shows that its
\(\bar\partial\)-data are rational and that the pole-type solution
of the corresponding replacement equation is rational. The algebraic
construction of the next meromorphic section in Theorem~5.3 therefore
remains in \(\mathbb C(z,\bar z)\), inductively at every actual
replacement level.

On the graph-regular set the highest Frenet coefficient is nonzero,
so dividing by it preserves \(\mathbb C(z,\bar z)\). The monic
vectors, graph matrices, and all Pfaffian expressions are therefore
rational in these variables. By Theorem~\ref{thm:psi-holomorphic}, their
\(\bar\partial\)-derivatives vanish. The kernel of
\(\partial_w\) in \(\mathbb C(z,w)\) is \(\mathbb C(z)\), proving
the first assertion. An identity of such polarized rational
functions holding on a dense part of the real diagonal is an
identity in the rational function field, proving the second.
\end{proof}

For each \(j\), clear the denominators of \(\boldsymbol\Psi_j\) and
divide out the common polynomial factor of its components. Write
\begin{equation}\label{eq:independent-saturation}
 \widehat{\boldsymbol\Psi}_j=\chi_j(z)\boldsymbol\Psi_j,\qquad
 Q_j=\norm{\widehat{\boldsymbol\Psi}_j}^2,\qquad
 D_j=\deg\widehat{\boldsymbol\Psi}_j.
\end{equation}
The resulting polynomial vector is primitive and has no common
zero on the affine chart. By the definition of its degree,
\(\zeta^{D_j}\widehat{\boldsymbol\Psi}_j(1/\zeta)\) is nonzero at
\(\zeta=0\), so it is a basepoint-free lift at infinity as well.
Set \(\chi_{-1}=1\), \(Q_{-1}=1\), and \(D_{-1}=0\).
Since \(\boldsymbol\Psi_0=f_0\) is already primitive, take \(\chi_0=1\).
Here \(D_0\) is the seed's actual degree. Each multiplier \(\chi_j\)
is chosen independently; its factors are retained when returning to
the raw equations.

\begin{proposition}[Saturated curvature equals actual energy]
\label{prop:saturated-energy}
For \(0\leq j\leq\ell\), the smooth curvature forms of the
saturated augmented curves satisfy
\begin{equation}\label{eq:saturated-energy}
 \partial\bar\partial\log Q_j=\sigma_j,\qquad
 \int_{\mathbb C}\sigma_j\,dx\,dy=\pi D_j.
\end{equation}
The function \(Q_j/(1+|z|^2)^{D_j}\) extends smoothly and
positively to \(\mathbb S^2\).
\end{proposition}

\begin{proof}
On the regular set, multiplication by the nonzero holomorphic
scalar \(\chi_j\) has zero mixed logarithmic derivative. Thus
\(\partial\bar\partial\log Q_j
=\partial\bar\partial\log P_j=\sigma_j\).
The left side is the smooth curvature of a basepoint-free
holomorphic curve. The right side is the smooth actual energy,
as established after Lemma~\ref{lem:rigidity-connection}.
Equality on the dense regular set therefore gives equality
	everywhere.

	To compute the degree integral directly, write
	\[
		\widehat{\boldsymbol\Psi}_j(z)
		=a_0+\cdots+a_{D_j}z^{D_j},
		\qquad a_{D_j}\ne0.
	\]
	At infinity,
	\[
		Q_j(z)
		=\norm{a_{D_j}}^2|z|^{2D_j}+O(|z|^{2D_j-1}),
	\]
	and hence
	\[
		\partial_r\log Q_j
		=\frac{2D_j}{r}+O(r^{-2}).
	\]
	Green's formula and
	\(\partial\bar\partial=\frac14\Delta\) give
	\[
		\int_{\mathbb C}\sigma_j\,dx\,dy
		=\lim_{R\to\infty}\frac14
		\int_{|z|=R}\partial_r\log Q_j\,ds
		=\pi D_j.
	\]
	The transition of the primitive polynomial lift at infinity also
	shows that its line bundle has degree \(-D_j\). The quotient
	\(Q_j/(1+|z|^2)^{D_j}\) extends smoothly and positively there.
\end{proof}

The factors \(\chi_j\) are independent. In particular,
\begin{equation}\label{eq:saturation-gauges}
 P_j=\frac{Q_j}{\chi_j\chi_j^*},\qquad
 B_j=\frac{H(\widehat{\boldsymbol\Psi}_{j+1},
                  \widehat{\boldsymbol\Psi}_{j-1})}
                  {\chi_{j+1}\chi_{j-1}^*}.
\end{equation}
These factors enter every conversion between the saturated data and
the raw recurrences. Their valuations at the round divisor are zero,
as used in \S\ref{sec:rigidity}.

\section{Effective Divisors and Growth of the Actual Degrees}
\label{sec:degree-growth}

Finite isotropy gives a strictly positive lower bound for the
increments of the actual degrees \(D_j\). The key input is the
holomorphic bilinear square on the middle Frenet line.

\begin{lemma}[The middle bilinear square]\label{lem:middle-square}
The restriction of \(b\) gives a nonzero holomorphic section
\[
 q=b|_{\mathcal L_{\ell+s+1}^{\otimes2}}
 \in H^0\bigl(\mathbb S^2,
             (\mathcal L_{\ell+s+1}^*)^{\otimes2}\bigr).
\]
If \(Z=\operatorname{div}(q)\), then
\begin{align}
 \deg\mathcal L_{\ell+s+1}&=-\tfrac12\deg Z\leq0,
 \label{eq:middle-square-degree}\\
 -\deg\mathcal L_j
 &=\tfrac12\deg Z+2(\ell+s+1-j)
   +\sum_{i=j}^{\ell+s}\epsilon_i,
 &&0\leq j\leq\ell+s.\label{eq:negative-seed-degrees}
\end{align}
\end{lemma}

\begin{proof}
In a smooth nonzero frame of the saturated line, write
\(\bar\partial e_{\ell+s+1}=\alpha e_{\ell+s+1}+v\) with
\(v\in\mathcal L_{\ell+s}\). By
\eqref{eq:short-mixed-orthogonality},
\(b(\mathcal L_{\ell+s},\mathcal L_{\ell+s+1})=0\). Hence
\[
 \bar\partial b(e_{\ell+s+1},e_{\ell+s+1})
 =2\alpha b(e_{\ell+s+1},e_{\ell+s+1}),
\]
which is exactly the induced holomorphic-section equation for
the dual square. This is a smooth saturated-bundle calculation,
so it also holds at seed ramification points.

To see that the section is nonzero, work on the monic regular set
and put \(F=b(f_0,f_{2(\ell+s)+2})\not\equiv0\), as required by the
finite seed class. Set \(p=2(\ell+s)+2\). For \(0\leq i<p\), the
short isotropy identities give
\(b(f_i,f_{p-i-1})=0\). Differentiating, using the Frenet equations,
and deleting the self-connection terms through the lower-order
vanishing pairings yields
\[
	b(f_{i+1},f_{p-i-1})+b(f_i,f_{p-i})=0.
\]
Thus induction from \(b(f_0,f_p)=F\) gives
\[
 b(f_i,f_{2(\ell+s)+2-i})=(-1)^iF,
 \qquad0\leq i\leq\ell+s+1.
\]
In particular
\(b(f_{\ell+s+1},f_{\ell+s+1})\not\equiv0\).
The divisor of this global nonzero holomorphic section has degree
\(-2\deg\mathcal L_{\ell+s+1}\), proving
\eqref{eq:middle-square-degree}. Summing
\eqref{eq:ramified-frenet-degrees} from \(j\) to \(\ell+s\) proves
\eqref{eq:negative-seed-degrees}, with every ramification term kept.
\end{proof}

The next step compares a retained line with a seed line, rather than
identifying them. That distinction permits the projection to vanish.

\begin{lemma}[Holomorphic highest-row projection]
\label{lem:highest-row-projection}
For \(0\leq j\leq\ell\), the smooth orthogonal projection restricts
to a nonzero holomorphic homomorphism for the induced structures,
\[
 \pi_j=\pi_{\mathcal L_j}|_{\underline v_j}:
 \underline v_j\longrightarrow\mathcal L_j.
\]
Its effective zero divisor \(\mathcal E_j\) satisfies
\begin{equation}\label{eq:projection-divisor-degree}
 \deg\underline v_j=\deg\mathcal L_j-\deg\mathcal E_j.
\end{equation}
\end{lemma}

The subtle point is that an arbitrary smooth orthogonal projection need
not be holomorphic. It becomes holomorphic here because the monic lift
\(v_j\) and the Frenet lift \(f_j\) are frames for the same induced
holomorphic structures, and \(\pi_jv_j=f_j\) on a dense set. The proof
below formulates precisely this annihilation of the \((0,1)\)-part.

\begin{proof}
The actual retained line, the saturated Frenet line, and the
orthogonal projection are globally smooth. On the regular set,
\eqref{eq:v-expansion} and the mixed pairing bound show precisely
that \(\pi_jv_j=f_j\): all lower positive rows and all conjugate
rows of \(v_j\) are orthogonal to \(\mathcal L_j\). The largest
mixed sum needed is \(2j-1\).

The graph connection and orthogonality imply
\(H(\bar\partial v_j,v_j)=0\). Thus \(v_j\) is a holomorphic
frame for the induced structure of its own line on this regular
set. The monic Frenet equation similarly makes \(f_j\) an
induced-holomorphic frame of \(\mathcal L_j\). Consequently,
\[
 (\bar\partial\pi_j)v_j
 =\bar\partial^{\mathcal L_j}(\pi_jv_j)
       -\pi_j(\bar\partial^{\underline v_j}v_j)=0.
\]
Since \(\bar\partial\pi_j\) is a globally smooth tensor and
vanishes on a dense set, it vanishes everywhere. The noncancelling
monic normalization shows that \(\pi_j\) is not the zero map.
The divisor-degree formula for a nonzero homomorphism of
holomorphic line bundles now gives
\eqref{eq:projection-divisor-degree}. We have neither excluded
the zeros of \(\pi_j\) in advance nor claimed that the monic
lifts are nonzero at every point.
\end{proof}

\begin{proposition}[Strict growth with all divisors retained]
\label{prop:saturated-degree-growth}
The actual saturated degrees satisfy
\begin{equation}\label{eq:actual-degree-increment}
 \boxed{\begin{aligned}
 D_j-D_{j-1}
 =\tfrac12\deg Z+2(\ell+s+1-j)
   +\sum_{i=j}^{\ell+s}\epsilon_i+\deg\mathcal E_j\\
 \geq2(\ell+s+1-j)\geq2,\qquad0\leq j\leq\ell.
 \end{aligned}}
\end{equation}
In particular \(D_{j+1}>D_{j-1}\) for \(0\leq j<\ell\).
\end{proposition}

\begin{proof}
On the regular set, the induced-holomorphic frame \(v_j\) has
Chern curvature represented, in the normalization used in
Proposition~\ref{prop:saturated-energy}, by
\(-\partial\bar\partial\log\norm{v_j}^2
=-(\sigma_j-\sigma_{j-1})\). Both sides extend as smooth forms
on the actual line bundle. Thus a degree is obtained by integrating
the curvature form and dividing by \(\pi\). Applying
Proposition~\ref{prop:saturated-energy} term by term, and using the
convention that the seed tautological line has negative degree,
gives
\begin{equation}\label{eq:retained-line-degree}
 \deg\underline v_j
 =-\frac1\pi\int_{\mathbb C}(\sigma_j-\sigma_{j-1})\,dx\,dy
 =D_{j-1}-D_j.
\end{equation}
Combine this with \eqref{eq:projection-divisor-degree} and
\eqref{eq:negative-seed-degrees}. Since \(j\leq\ell\), every
term discarded in the lower bound is nonnegative and
\(\ell+s+1-j\geq1\).
\end{proof}

Frenet ramification, zeros of the highest-row projection, and zeros of
the middle square contribute nonnegative terms to this lower bound.

\section{Constant Curvature and Raw Backward Rigidity}
\label{sec:rigidity}

An immersed constant-curvature metric on \(\mathbb S^2\) has \(K>0\)
and admits a spherical affine coordinate in which
\(g_\ell=4K^{-1}(1+|z|^2)^{-2}\,dz\,d\bar z\).
We use this coordinate throughout the rigidity argument.

The proof below has three layers. First, the constant-curvature metric
rounds only the two terminal norms \(Q_{\ell-1},Q_\ell\). Second,
the exact raw norm and source equations are converted, through
valuations along \(1+zw=0\), into integer inequalities. Third,
those inequalities propagate a nonzero deviation from the top of the
tower down to \(B_0=0\), where the scalar bottleneck gives a
contradiction. Only after this backward argument is complete do we
identify the seed as a Veronese curve.

\begin{lemma}[Terminal rounding]\label{lem:rigidity-rounding}
In this common coordinate,
\begin{equation}\label{eq:terminal-saturated-rounding}
 \begin{gathered}
 \frac2K=D_\ell+D_{\ell-1},\\
 Q_j=C_j(1+|z|^2)^{D_j},\qquad
 C_j>0,\quad j=\ell-1,\ell.
 \end{gathered}
\end{equation}
\end{lemma}

\begin{proof}
The actual metric identity \eqref{eq:rigidity-metric} gives
\(\sigma_\ell+\sigma_{\ell-1}
=(2/K)(1+|z|^2)^{-2}\).
Integrating and using \eqref{eq:saturated-energy} proves the
degree-curvature identity. The smooth function
\[
 \log\frac{Q_\ell Q_{\ell-1}}
              {(1+|z|^2)^{D_\ell+D_{\ell-1}}}
\]
is harmonic on the sphere, so it is constant. Polarize this polynomial
identity by replacing \(\bar z\) with \(w\), obtaining
\[
	Q_\ell(z,w)Q_{\ell-1}(z,w)
	=C(1+zw)^{D_\ell+D_{\ell-1}},
	\qquad C>0.
\]
Write
\[
	\widehat{\boldsymbol\Psi}_j(z)
	=\sum_{a=0}^{D_j}u_{j,a}z^a,
	\qquad u_{j,D_j}\ne0.
\]
The polarized norm has a nonzero
\(z^{D_j}w^{D_j}\)-coefficient, namely
\(\norm{u_{j,D_j}}^2\). Hence \(Q_j(z,w)\) has bidegree exactly
\((D_j,D_j)\). The polynomial \(1+zw\) is irreducible in
\(\mathbb C[z,w]\); for example its quotient ring is
\(\mathbb C[z,z^{-1}]\). Unique factorization therefore writes
each factor as a constant times a nonnegative integral power of
\(1+zw\). Comparison of the two \(z\)-degrees fixes the exponents
to be \(D_\ell\) and \(D_{\ell-1}\). Restricting back to
\(w=\bar z\) gives the two round powers in
\eqref{eq:terminal-saturated-rounding}. Their constants are positive
on the real diagonal.
\end{proof}

To propagate this roundness through the exact raw equations, we
retain the saturation factors and measure orders along the round
divisor. In \(\mathbb C(z,w)\), write
\(\operatorname{ord}_{1+zw}\) for the order along \(1+zw=0\),
with \(\operatorname{ord}_{1+zw}(0)=+\infty\). In formulas with polarized
data, \(\partial=\partial_z\) and \(\bar\partial=\partial_w\).
The involution \(F^*(z,w)=\overline{F(\bar w,\bar z)}\) fixes
\(1+zw\), so
\(\operatorname{ord}_{1+zw}(F^*)=\operatorname{ord}_{1+zw}(F)\).
Every nonzero one-variable
rational function has valuation zero. In particular the factors
\(\chi_j\) and \(\chi_j^*\) in
\eqref{eq:saturation-gauges} do not change the valuation.

All orders below are taken in the local ring of the irreducible
divisor \(1+zw=0\). Thus a denominator which is nonzero there is a
unit and contributes order zero, while orders add under multiplication.
This is the precise sense in which polynomial divisibility controls
the poles of all logarithmic derivatives appearing in the raw
equations. The valuation calculation is insensitive to common
one-variable factors. Although the saturation factors come from
different levels, each is a rational function of one variable and
therefore has order zero along \(1+zw=0\). Independence of the
\(\chi_j\) matters when translating back to the \(Q_j\), but not
when evaluating their orders.

Put
\begin{equation}\label{eq:raw-valuations}
 t_i=\operatorname{ord}_{1+zw}(P_i),\quad t_{-1}=0,\qquad
 \upsilon_j=\operatorname{ord}_{1+zw}(B_j)\quad\hbox{when }B_j\ne0.
\end{equation}

\begin{lemma}[Raw valuation bounds and derivative orders]
\label{lem:raw-valuation-calculus}
For nonzero mixed kernels,
\begin{equation}\label{eq:raw-valuation-bounds}
 0\leq t_i\leq D_i,\qquad
 0\leq\upsilon_j\leq D_{j-1},\qquad
 B_0\ne0\Longrightarrow\upsilon_0=0.
\end{equation}
For nonzero rational functions \(P,B\), put
\(t=\operatorname{ord}_{1+zw}(P)\) and
\(\upsilon=\operatorname{ord}_{1+zw}(B)\). Then
\begin{align}
 \operatorname{ord}_{1+zw}(P^2\partial\bar\partial\log P)&=2t-2
 &&(t>0),
 \label{eq:raw-curvature-order}\\
 \operatorname{ord}_{1+zw}(-P\bar\partial B+(\bar\partial P)B)
 &=t+\upsilon-1
 &&(t\ne\upsilon),
 \label{eq:raw-source-exact-order}\\
 \operatorname{ord}_{1+zw}(P\partial B-(\partial P)B)
 &\geq t+\upsilon-1.
 \label{eq:raw-source-bound}
\end{align}
\end{lemma}

\begin{proof}
The polarized \(Q_i\) has bidegree \((D_i,D_i)\), and the
numerator of \(B_j\) in \eqref{eq:saturation-gauges} has
bidegree at most \((D_{j+1},D_{j-1})\). The denominators have
valuation zero. Polynomial divisibility therefore gives the first
two bounds. The bottom kernel is a rational function of \(z\)
alone, so a nonzero \(B_0\) has valuation zero.

Write \(P=(1+zw)^tp\) and \(B=(1+zw)^\upsilon b\),
where \(p,b\) are units along the divisor. Their derivatives
have no poles there, whereas
\(\partial\bar\partial\log(1+zw)=(1+zw)^{-2}\).
For \(t>0\), the leading term \(t(1+zw)^{-2}\) proves
\eqref{eq:raw-curvature-order}. For the source equation,
\[
 -P\bar\partial B+(\bar\partial P)B
 =(1+zw)^{t+\upsilon-1}
 \bigl((t-\upsilon)zpb
       +(1+zw)((\bar\partial p)b-p\bar\partial b)\bigr).
\]
Since \(z\) is a unit along the divisor, the leading coefficient
cannot vanish when \(t\ne\upsilon\). The same expansion
with \(\partial\) gives the general lower bound
\eqref{eq:raw-source-bound}.
\end{proof}

We can now rule out a nonzero deviation between two round levels.
The source equation carries its nonvanishing backwards until it
contradicts the bottom kernel bound.

\begin{lemma}[Raw backward rigidity]\label{lem:raw-backward-rigidity}
Suppose two adjacent saturated norms are round:
\(Q_k=C_k(1+|z|^2)^{D_k}\) and
\(Q_{k+1}=C_{k+1}(1+|z|^2)^{D_{k+1}}\),
where \(0\leq k<\ell\). Then \(B_k=0\). If \(k\geq1\),
\begin{equation}\label{eq:backward-rounding-step}
 D_{k-1}=2D_k-D_{k+1}-2,\qquad
 Q_{k-1}=C_{k-1}(1+|z|^2)^{D_{k-1}},\quad C_{k-1}>0.
\end{equation}
For \(k=0\) the degree identity is \(D_1=2D_0-2\).
\end{lemma}

The proof should be read in three stages. First, if
\(B_k\ne0\), the strict degree growth forces its valuation to be
strictly smaller than the curvature valuation in the norm equation;
this produces an exact initial descent relation. Second, the source
equation carries nonvanishing from level \(j\) to level \(j-1\), while
the norm equation controls how quickly the relevant valuations can
decrease. Finally, at \(j=1\), the scalar nature of \(B_0\) makes
the bottom kernel valuation zero, contradicting the descent.

\begin{proof}
Rounding and the zero valuations of the saturation factors give
\(t_k=D_k\), \(t_{k+1}=D_{k+1}\), and
\(\partial\bar\partial\log P_k=D_k(1+zw)^{-2}\).
Suppose \(B_k\ne0\). By
\eqref{eq:actual-degree-increment} and
\eqref{eq:raw-valuation-bounds},
\[
 \upsilon_k\leq D_{k-1}\leq D_k-2=t_k-2.
\]
The kernel term in \eqref{eq:rigidity-N} then has strictly lower
valuation than the curvature term, so they cannot cancel. Hence
\begin{equation}\label{eq:initial-descent-order}
 t_{k-1}+D_{k+1}=2\upsilon_k.
\end{equation}
If \(k=0\), the bottom bound \(\upsilon_0=0\) forces
\(D_1=0\), contradicting strict degree growth. For \(k\geq1\),
\(D_{k+1}\geq D_k+2\) and the preceding inequality give
\(0\leq t_{k-1}<t_k<t_{k+1}\).

We now propagate the following conditions backwards:
\begin{equation}\label{eq:descent-invariants}
 B_j\ne0,\qquad\upsilon_j\leq t_j-2,\qquad
 0\leq t_{j-1}<t_j\leq t_{j+1},\qquad1\leq j\leq k.
\end{equation}
They hold at \(j=k\). The norm equation and
\eqref{eq:raw-curvature-order} give
\(2\upsilon_j=t_{j-1}+t_{j+1}\), whence
\(\upsilon_j>t_{j-1}\). Thus the right side of the source
equation \eqref{eq:rigidity-C} has the exact finite valuation
\(t_{j-1}+\upsilon_j-1\), by
\eqref{eq:raw-source-exact-order}. It is nonzero, so
\(B_{j-1}\ne0\). The lower bound on the left side gives
\begin{equation}\label{eq:strict-valuation-descent}
 \upsilon_{j-1}\leq t_{j-1}+\upsilon_j-t_j
                         \leq t_{j-1}-2.
\end{equation}
In particular \(t_{j-1}\geq2\). For \(j\geq2\), applying
the norm equation at \(j-1\) gives
\[
 t_{j-2}+t_j=2\upsilon_{j-1}\leq2t_{j-1}-4,
\]
and therefore
\[
 t_{j-2}\leq 2t_{j-1}-4-t_j<t_{j-1}-4<t_{j-1}<t_j.
\]
This proves \eqref{eq:descent-invariants} at the next index.
At \(j=1\), \eqref{eq:strict-valuation-descent} instead gives
\(B_0\ne0\), \(t_0\geq2\), and \(t_0<t_1\).
But \(\upsilon_0=0\) and the dominating bottom kernel term in
\eqref{eq:rigidity-N} imply
\(t_{-1}+t_1=2\upsilon_0=0\), a contradiction.
Therefore \(B_k=0\).

The exact norm equation now reads
\[
 P_{k-1}=\frac{P_k^2}{P_{k+1}}\,D_k(1+zw)^{-2}.
\]
Taking mixed logarithmic derivatives removes all one-variable
gauge factors and gives
\[
 \partial\bar\partial\log P_{k-1}
 =(2D_k-D_{k+1}-2)(1+zw)^{-2}.
\]
For \(k\geq1\), this is an equality of the smooth saturated
curvature forms. Integrating proves the degree identity in
\eqref{eq:backward-rounding-step}; the smooth harmonic quotient
\(\log(Q_{k-1}/(1+|z|^2)^{D_{k-1}})\) is then constant. For \(k=0\),
the left side is zero because \(P_{-1}=1\), giving
\(D_1=2D_0-2\).
\end{proof}

\begin{theorem}[Reverse rigidity in the raw gauge]
\label{thm:rigidity-reverse}
For a constant-curvature target, all kernels \(B_j\) and all
deviations \(x_j\) vanish. In the same spherical coordinate,
every saturated norm is round and
\begin{equation}\label{eq:actual-degree-recurrence}
 Q_j=C_j(1+|z|^2)^{D_j},\qquad
 D_{j+1}=2D_j-D_{j-1}-2,\qquad D_{-1}=0.
\end{equation}
On the regular set \(v_j=f_j\); the corresponding saturated
lines agree everywhere.
\end{theorem}

\begin{proof}
Begin with Lemma~\ref{lem:rigidity-rounding} and apply
Lemma~\ref{lem:raw-backward-rigidity} for
\(k=\ell-1,\ell-2,\ldots,0\). The degree-growth hypothesis
needed at every step was proved before using constant curvature.
Thus every \(B_j\) vanishes and all \(Q_j\) round. Since
\(B_j=x_jP_j\) and \(P_j>0\) on the regular set, every
\(x_j\) vanishes there. Proposition~\ref{prop:faithful-frenet}
gives \(v_j=f_j\) on that set. Equality of smooth saturated
lines extends everywhere by continuity of their projections.
\end{proof}

At depth one the two terminal norms are already \(Q_0,Q_1\);
the bottom step forces the scalar grading \(B_0\) to vanish.
At depth two the metric initially rounds only \(Q_1,Q_2\).
The first descent forces \(B_1=0\) and rounds \(Q_0\), after which
the bottom step forces \(B_0=0\). The same descent applies at every
depth.

\section{Veronese Representation and Automatic Unramifiedness}
\label{sec:automatic-unramifiedness}

The faithful-tower criterion and reverse rigidity now give the
Veronese representation directly. Only the seed's projective span
is needed to identify its degree.

The proof has three steps: coefficient comparison for the rounded
seed norm identifies its degree and embedding; the Frenet equations
then identify every retained lift with a Veronese lift; finally,
the rank-two target's adjacent Gauss layers are read off from the
positive and conjugate Veronese rows. This last step is where the
otherwise abstract factorization turns into total unramifiedness
of the original map.

\begin{theorem}[Veronese representation in the projective span]
\label{thm:veronese-representation}
Under the hypotheses of Theorem~\ref{thm:rigidity-reverse},
\(D_0=m\) and there is a complex linear isometric embedding
\(U:\mathbb C^{m+1}\hookrightarrow\mathbb C^n\) such that
\begin{equation}\label{eq:veronese-isometric-embedding}
 f_0=\sqrt{C_0}\,UV_0^{(m)},\qquad U^*U=I_{m+1}.
\end{equation}
Every seed Frenet homomorphism is unramified. The monic retained
lifts are \(v_j=f_j=\sqrt{C_0}\,UV_j^{(m)}\), \(0\leq j\leq\ell\),
and
\begin{equation}\label{eq:veronese-real-layer}
 \underline\varphi
 =\operatorname{span}(UV_\ell^{(m)})\oplus
   \operatorname{span}(\overline{UV_\ell^{(m)}}).
\end{equation}
With the saturated lifts chosen in the raw normalization, the
degrees and norms satisfy
\begin{equation}\label{eq:rigidity-all-round}
 D_j=d_j=(j+1)(m-j),\qquad
 P_j=C_j(1+|z|^2)^{d_j},\qquad0\leq j\leq\ell.
\end{equation}
\end{theorem}

\begin{proof}
By Proposition~\ref{prop:faithful-frenet} and
Theorem~\ref{thm:rigidity-reverse}, \(v_j=f_j\) and
\(\norm{f_0}^2=C_0(1+|z|^2)^{D_0}\). Writing
\(f_0(z)=\sum_{h=0}^{D_0}a_hz^h\) and comparing coefficients gives
\begin{equation}\label{eq:seed-coefficient-orthogonality}
 H(a_h,a_k)=C_0\binom{D_0}{h}\delta_{hk}.
\end{equation}
Indeed, expanding
\(\norm{f_0(z)}^2=C_0(1+|z|^2)^{D_0}\) gives the diagonal binomial
coefficients and the vanishing of all mixed terms. In particular,
no coefficient vector \(a_h\) is zero.

The projective span of the seed is, by definition, the linear span
of its coefficient vectors and has dimension \(m+1\). Since all
\(D_0+1\) coefficients are nonzero and mutually orthogonal, they are
linearly independent and span a space of dimension \(D_0+1\).
Consequently \(D_0=m\). The normalization
\[
	U(e_h)=\frac{a_h}{\sqrt{C_0\binom{D_0}{h}}},
	\qquad0\leq h\leq D_0,
\]
is an isometric embedding \(U:\mathbb C^{m+1}\hookrightarrow
\mathbb C^n\) and gives
\eqref{eq:veronese-isometric-embedding}. Frenet differentiation
commutes with \(U\), and the standard Veronese formulas give
\begin{equation}\label{eq:faithful-veronese-densities}
 \lambda_i=\frac{(i+1)(m-i)}{(1+|z|^2)^2},\qquad0\leq i<m.
\end{equation}
These intrinsic densities are positive on the sphere, proving total
unramifiedness. The retained-line formula and
\eqref{eq:rigidity-all-round} follow from \(v_j=f_j\),
the degree recurrence \eqref{eq:actual-degree-recurrence}, and
the same Veronese norm formulas. The resulting raw osculating
lifts are primitive, so the saturation factors can be chosen as one.
\end{proof}

The possibly rectangular embedding \(U\) records the complex seed
span; real fullness will be treated in Section~\ref{sec:real-gram}.
To transfer the positive Veronese densities to the rank-two target,
we identify its actual adjacent Gauss bundles.

\begin{proposition}[Automatic target unramifiedness]
\label{prop:target-automatic-unramifiedness}
After the identification in Theorem~\ref{thm:veronese-representation},
the actual positive Gauss bundles are
\begin{align}
 E_j&=\mathcal L_{\ell+j}\oplus
                 \overline{\mathcal L_{\ell-j}},
 &&0\leq j\leq\ell,\label{eq:faithful-gauss-middle}\\
 E_{\ell+h}&=\mathcal L_{2\ell+h},
 &&1\leq h\leq m-2\ell.\label{eq:faithful-gauss-tail}
\end{align}
The target's two-sided harmonic sequence is totally unramified.
\end{proposition}

\begin{proof}
Start with \(E_0=\mathcal L_\ell\oplus\overline{\mathcal L_\ell}\)
and differentiate. In each middle layer, the positive row moves
to \(\mathcal L_{\ell+j+1}\) and the conjugate row moves to
\(\overline{\mathcal L_{\ell-j-1}}\). Their orthogonality to
the current layer and to one another uses mixed index sums at
most \(2\ell+1\), which lie in
\eqref{eq:short-mixed-orthogonality}. Their self-connection
parts are removed by the current orthogonal projection.

For \(j<\ell\), in the notation
\(E_j=\mathcal L_{\ell+j}\oplus
\overline{\mathcal L_{\ell-j}}\), the two normal columns are,
up to nonzero scalar factors,
\[
	f_{\ell+j+1},
	\qquad
	\lambda_{\ell-j-1}\overline{f_{\ell-j-1}}.
\]
They are orthogonal because
\(H(f_{\ell+j+1},\overline{f_{\ell-j-1}})
=b(f_{\ell+j+1},f_{\ell-j-1})\), whose index sum is \(2\ell\).
These are the two columns whose normalized Gram determinant is
computed next. At \(j=\ell\), the conjugate row contains the factor
\(\lambda_{-1}=0\), while the remaining row is \(f_{2\ell+1}\).
This is the explicit rank-two-to-rank-one interface. The subsequent
sequence is the ordinary projective tail. This proves the formulas
on a dense set, and equality of saturated bundles gives them globally.

For \(j<\ell\), in unitary frames the two normal columns have
Gram determinant density
\[
 \lambda_{\ell+j}\lambda_{\ell-j-1}>0.
\]
At \(j=\ell\), the sole nonzero normal column has squared norm
\(\lambda_{2\ell}>0\). All later nonterminal maps have the
positive density \(\lambda_i\) of the projective tail.
At infinity, use the smooth local coordinate \(\zeta=1/z\). The
transition formula \eqref{eq:transition-frenet-v} expresses the
same Frenet line in a nonzero local frame, and changing a frame
cannot change the rank of a homomorphism. Thus the positivity is
not restricted to the affine coordinate and the same conclusion
holds in a neighbourhood of infinity.
The admissible range ensures that these indices exist. Thus
every positive-direction homomorphism has its generic rank
everywhere. Reality gives \(E_{-j}=\overline{E_j}\) and the
same assertion in the negative direction. This is exactly total
unramifiedness in Definition~\ref{def:unramified}.
\end{proof}

Combining
\eqref{eq:terminal-saturated-rounding} with the now proved standard
degrees gives
\[
 \frac2K=D_\ell+D_{\ell-1}
 =(\ell+1)(m-\ell)+\ell(m-\ell+1)
 =(2\ell+1)m-2\ell^2,
\]
as claimed. At depths one and two this is
\(K=2/(3m-2)\) and \(K=2/(5m-8)\), respectively.

\section{Real Gram Data, Converse, and Examples}\label{sec:real-gram}

The real extrinsic data of the embedding are encoded by the symmetric
Gram matrix \(\mathsf R=U^{\transp}U\). Following
\citep{chi2021structure}, we express finite isotropy and real fullness
as algebraic conditions on this matrix. The matrix retains the
extrinsic data beyond the intrinsic Veronese description.

The converse proof follows the reverse order. The algebraic kernel
criterion first converts finite isotropy into divisibility of a
polynomial kernel. The real Gram matrix then controls realizability
and real fullness. We next rebuild the Veronese real layers, check
that their forward and backward replacements do not cancel, and use
BEW uniqueness to identify the depth. The final paragraph converts
the positive Veronese metric back to the stated curvature.

    \begin{proposition}[Real Gram data]
    	\label{prop:real-gram}
    	Fix the standard Veronese basis, and let
    	\(U:\mathbb C^{m+1}\hookrightarrow\mathbb C^n\) be the embedding
    	in Theorem~\ref{thm:veronese-representation}. Put
    	\begin{equation}
    		\mathsf R=U^{\transp}U
    		\in\mathbb{C}^{(m+1)\times(m+1)}.
    	\end{equation}
    	Then \(\mathsf R=\mathsf R^{\transp}\) and
    	\(\mathsf R^*\mathsf R\leq I_{m+1}\). For all
    	\(u,w\in\mathbb{C}^{m+1}\),
    	\begin{equation}\label{eq:rigidity-realgram}
    		b(Uu,Uw)=u^{\transp}\mathsf Rw.
    	\end{equation}
    	The initial finite isotropy conditions are exactly
    	\begin{equation}\label{eq:realgram-isotropy}
    	 \begin{aligned}
    	 &(V_0^{(m)})^{\transp}\mathsf RV_i^{(m)}=0,
    	 &&0\leq i\leq2(\ell+s)+1,\\
    	 &(V_0^{(m)})^{\transp}\mathsf RV_{2(\ell+s)+2}^{(m)}
    	 \not\equiv0.
    	 \end{aligned}
    	\end{equation}
    	Equivalently, with \(p=2(\ell+s)+2\) and
    	\begin{equation}\label{eq:realgram-kernel}
    	 \mathcal K_{\mathsf R}(z,t)
    	 =(V_0^{(m)}(z))^{\transp}\mathsf RV_0^{(m)}(t),
    	\end{equation}
    	they are the finite algebraic conditions
    	\begin{equation}\label{eq:realgram-divisibility}
    	 (t-z)^p\mid\mathcal K_{\mathsf R}(z,t),
    	 \qquad
    	 (t-z)^{p+1}\nmid\mathcal K_{\mathsf R}(z,t).
    	\end{equation}
    	Put
    	\begin{equation}\label{eq:realgram-real-matrix}
    	 \mathcal G(\mathsf R)=\frac12
    	 \begin{pmatrix}
    	 I+\operatorname{Re}\mathsf R&\operatorname{Im}\mathsf R\\
    	 \operatorname{Im}\mathsf R&I-\operatorname{Re}\mathsf R
    	 \end{pmatrix}.
    	\end{equation}
    	Then the target real layer is linearly full precisely when
    	\begin{equation}\label{eq:realgram-fullness}
    	 \operatorname{rank}\mathcal G(\mathsf R)=n.
    	\end{equation}
    	More generally, any symmetric contraction \(\mathsf R\) is
    	realized by a complex isometric embedding into \(\mathbb C^n\)
    	if and only if \(\operatorname{rank}\mathcal G(\mathsf R)\leq n\).
    	Two complex isometric embeddings \(U,U'\), with the same
    	domain basis, satisfy \(U'=OU\) for some \(O\in O(n)\) if and
    	only if \(U^{\transp}U=U'^{\transp}U'\).
    \end{proposition}

    \begin{proof}
    	Symmetry and \eqref{eq:rigidity-realgram} follow directly from
    	\(b(u,w)=u^{\transp}w\). Since
    	\(\mathsf R=(\overline U)^*U\) and both \(U\) and
    	\(\overline U\) are isometric embeddings,
    	\(\mathsf R\) is a contraction. Substituting
    	\(f_i=\sqrt{C_0}\,UV_i^{(m)}\) into the defining conditions of
    	\(H_n^{\ell+s}\) proves \eqref{eq:realgram-isotropy}.

    	The Frenet vectors are triangular combinations of ordinary
    	derivatives, with leading coefficient one:
    	\[
    		f_i(z)
    		=\partial_z^if_0(z)
    		+\sum_{h=0}^{i-1}c_{i,h}(z)\partial_z^hf_0(z).
    	\]
    	Consequently \(b(f_0,f_i)\) equals
    	\(\partial_t^i\mathcal K_{\mathsf R}(z,t)|_{t=z}\) plus
    	lower derivatives. The vanishing for every \(0\leq i<p\) is
    	therefore equivalent, by induction on \(i\), to the vanishing
    	of the first \(p\) \(t\)-derivatives at \(t=z\). Taylor
    	expansion in \(t-z\) proves \eqref{eq:realgram-divisibility},
    	including the nonzero condition of order \(p\).

    	For the remaining assertions, write \(U=X+\mathrm{i}Y\), with real
    	\(n\times(m+1)\) matrices \(X,Y\). The identities
    	\(U^*U=I\) and \(U^{\transp}U=\mathsf R\) give
    	\begin{equation}\label{eq:rigidity-realblocks}
    		X^{\transp}X
    		=\frac{I+\operatorname{Re}\mathsf R}{2},
    		\quad
    		Y^{\transp}Y
    		=\frac{I-\operatorname{Re}\mathsf R}{2},
    		\quad
    		X^{\transp}Y
    		=\frac{\operatorname{Im}\mathsf R}{2}.
    	\end{equation}
    	Consequently,
    	\(\mathcal G(\mathsf R)=[X\ Y]^{\transp}[X\ Y]\).
    	The values of \(f_\ell\) span the same complex space as those of
    	\(f_0\). Indeed, the backward Frenet formula
    	\[
    		\bar\partial f_i
    		=-\lambda_{i-1}f_{i-1}
    		\qquad(1\leq i\leq\ell)
    	\]
    	and \(\lambda_{i-1}>0\) show that a constant covector
    	annihilating \(f_\ell\) annihilates \(f_{\ell-1}\), and then
    	all \(f_i\) down to \(f_0\). Conversely, the forward formula
    	\[
    		\partial f_i
    		=f_{i+1}
    		+(\partial\log\norm{f_i}^2)f_i
    	\]
    	shows by induction that a constant covector annihilating \(f_0\)
    	annihilates every \(f_i\). Hence the real layers span exactly
    	the real column span of \([X\ Y]\). Its dimension is
    	\(\operatorname{rank}\mathcal G(\mathsf R)\), proving
    	\eqref{eq:realgram-fullness}. In particular, full real
    	linearity requires
    	\[
    		n\leq2(m+1),
    	\]
    	because \([X\ Y]\) has \(2(m+1)\) columns.

    	For the realizability assertion, set
    	\(A=\operatorname{Re}\mathsf R\),
    	\(B=\operatorname{Im}\mathsf R\), and
    	\(T=\left(\begin{smallmatrix}A&B\\B&-A\end{smallmatrix}\right)\).
    	This is a real symmetric matrix representing
    	\(v\mapsto\mathsf R\overline v\). Since \(\mathsf R\) is
    	a contraction, the eigenvalues of \(T\) lie in \([-1,1]\),
    	so \(\mathcal G(\mathsf R)=(I+T)/2\) is positive
    	semidefinite. It has a real Gram factorization
    	\([X\ Y]^{\transp}[X\ Y]\) in \(\mathbb R^n\) precisely
    	when its rank is at most \(n\). Its blocks then give
    	\(U=X+\mathrm{i}Y\) with \(U^*U=I\) and
    	\(U^{\transp}U=\mathsf R\).

    	Finally, \(\mathsf R\) determines all real Gram products of the
    	columns of \(X\) and \(Y\). For two embeddings with the same
    	\(\mathsf R\), equality of these real Gram matrices defines an
    	isometry between the real spans of the corresponding columns.
    	It extends to an orthogonal transformation of \(\mathbb R^n\),
    	giving \(U'=OU\). The converse follows from
    	\(O^{\transp}O=I\).
    \end{proof}

    \begin{corollary}[Algebraic converse]\label{cor:veronese-converse}
    	Fix integers \(s\geq0\), \(\ell\geq1\), and
    	\(m\geq p:=2(\ell+s)+2\). Let \(\mathsf R\) be a symmetric
    	contraction satisfying \eqref{eq:realgram-divisibility} and
    	\eqref{eq:realgram-fullness}, and choose an isometric embedding
    	\(U:\mathbb C^{m+1}\hookrightarrow\mathbb C^n\) realizing it.
    	Then
    	\[
    	 \underline\varphi
    	 =\operatorname{span}(UV_\ell^{(m)})\oplus
    	 \operatorname{span}(\overline{UV_\ell^{(m)}})
    	\]
    	defines a linearly full, irreducible, totally
    	unramified conformal minimal immersion of isotropy order
    	\(2s+1\), replacement depth \(\ell\), and constant Gaussian
    	curvature \eqref{eq:main-curvature}. Together with
    	Theorem~\ref{thm:main-veronese-classification}, this describes
    	exactly the family in the stated range.
    \end{corollary}

    \begin{proof}
    	Set \(f_j=UV_j^{(m)}\) and
    	\(\mathcal L_j=\operatorname{span}(f_j)\). By the kernel
    	criterion, \(b(f_0,f_i)=0\) for \(i<p\) and
    	\(b(f_0,f_p)\not\equiv0\). Differentiation of these identities
    	and the Frenet equations give
    	\[
    	 b(f_a,f_b)=0\quad(a+b<p),\qquad
    	 b(f_a,f_{p-a})=(-1)^a b(f_0,f_p)
    	 \quad(0\leq a\leq p).
    	\]
    	In particular, the short windows needed for each current
    	projection and adjacent image are orthogonal; no distant-block
    	orthogonality is required. For
    	\(V_j=\mathcal L_j\oplus\overline{\mathcal L_j}\),
    	\(0\leq j\leq\ell\), the two normal columns are
    	\[
    	 A'_{V_j}(f_j)=f_{j+1},\qquad
    	 A'_{V_j}(\overline{f_j})=-\lambda_{j-1}\overline{f_{j-1}}
    	 \quad(j\geq1),
    	\]
    	where the second column is zero at \(j=0\).
    	The Frenet equations
    	show that these columns define a holomorphic second
    	fundamental form, so \(V_j\) is harmonic. They also give the
    	faithful replacements with
    	\(V_j=\mathcal L_j\oplus\overline{\mathcal L_j}\) and
    	\(U_j=\mathcal L_{j+1}\oplus\overline{\mathcal L_j}\),
    	starting at the harmonic real mixed pair \(V_0\).

    	The actual forward Gauss bundles of \(V_\ell\) are
    	\[
    	 \begin{aligned}
    	 G^{(k)}(V_\ell)
    	 &=\mathcal L_{\ell+k}\oplus
    	   \overline{\mathcal L_{\ell-k}},
    	 &&0\leq k\leq\ell,\\
    	 G^{(\ell+a)}(V_\ell)&=\mathcal L_{2\ell+a},
    	 &&1\leq a\leq m-2\ell.
    	 \end{aligned}
    	\]
    	Every nonterminal connection has full row rank because
    	\(\norm{f_{j+1}}^2/\norm{f_j}^2
    	=(j+1)(m-j)/(1+|z|^2)^2>0\); the same assertion holds at
    	infinity in a smooth local frame. The first nonorthogonal
    	Gauss bundle occurs at \(k=p-2\ell=2s+2\), by the displayed
    	bilinear identities. Hence the isotropy order is \(2s+1\).
    	The same calculation at each earlier real layer gives order
    	\(p-2j-1\) for \(V_j\). Consequently none of the forward
    	replacements \(V_j\mapsto U_j\) or backward replacements
    	\(U_j\mapsto V_{j+1}\), \(0\leq j<\ell\), cancels with its
    	inverse. Indeed, the first nonorthogonal transform of \(V_j\)
    	occurs only after \(p-2j\) steps.

    	We have therefore constructed, from the real mixed pair
    	\(V_0\), exactly a sequence
    	\[
    		V_0,\ U_0,\ V_1,\ldots,U_{\ell-1},\ V_\ell
    	\]
    	of harmonic forward and backward replacements satisfying
    	conditions (i)--(vi) of Theorem~4.7 of
    	\citep{bahy1989explicit}. Their isotropy orders decrease by
    	one at each half-step, from \(2(\ell+s)+1\) at \(V_0\) to
    	\(2s+1\) at \(V_\ell=\varphi\). Theorem~4.7 states that this
    	factorization is unique. Hence the canonical noncancelling
    	replacement depth of \(\varphi\) is exactly \(\ell\), and no
    	additional inverse pair has been inserted. The rank of the
    	first Gauss bundle is two, proving irreducibility, and
    	\eqref{eq:realgram-fullness} proves linear fullness.

    	Harmonicity on the sphere implies weak conformality, since
    	its holomorphic Hopf differential must vanish. The
    	Grassmannian metric is
    	\[
    	 g_\ell=\frac{2\bigl((2\ell+1)m-2\ell^2\bigr)}
    	 {(1+|z|^2)^2}\,\mathrm{d}z\,\mathrm{d}\overline z,
    	\]
    	which is positive everywhere, so there are no branch points.
    	The map is therefore a conformal minimal immersion with the
    	claimed constant curvature.
    \end{proof}

\subsection{Conclusion of the classification}\label{sec:main-veronese}

\begin{proof}[Proof of Theorem~\ref{thm:main-veronese-classification}]
The BEW construction supplies the actual tower and the seed class.
The short-window graph calculus yields the raw Pfaffian equations
and the actual smooth energies. Independent saturation gives
Proposition~\ref{prop:saturated-energy}, while finite isotropy and
the highest-row projection give
Proposition~\ref{prop:saturated-degree-growth}. Neither step
uses an unramifiedness assumption. Constant curvature and raw backward
rigidity then give
Theorem~\ref{thm:rigidity-reverse}. The seed coefficient comparison
in Theorem~\ref{thm:veronese-representation} proves the degree,
embedding, and seed unramifiedness; Proposition~\ref{prop:target-automatic-unramifiedness}
independently proves target unramifiedness. The degree-curvature
identity gives \eqref{eq:main-curvature}, and
Proposition~\ref{prop:real-gram} gives the necessary real Gram
constraints. Corollary~\ref{cor:veronese-converse} proves their
sufficiency. This completes both directions.
\end{proof}

\subsection{Comparison with low-dimensional classifications}
\label{subsec:known-examples}

The \(Q_2\) and \(Q_3\) models of
\citep{jiao2013q2,li2016q3} are recorded explicitly in
\citep[Corollaries~4.7--4.8]{jiao2017classification}; the latter
corollary reproduces Theorem~4.9 of \citep{li2016q3}. For \(Q_4\)
we use \citep[\S\S4--5]{jiao2021q4}, and for
\(Q_5\) we use \citep[Lemma~5 and Theorems~1--2]{jiao2019q5}.
The irreducible classifications quoted for \(Q_3,Q_4,Q_5\) assume
total unramifiedness. In the \(Q_5\), \(r=1\) case, identification
of the model in Theorem~1 also assumes polynomiality of the indicated
second Frenet norm; Theorem~2 determines its curvature without that
additional assumption.

For the table, write
\[
 \Phi_j^{(m)}=[UV_j^{(m)}]\oplus[\overline{UV_j^{(m)}}],
\]
where \([v]\) denotes the complex line spanned by \(v\), with the
admissible embeddings \(U\) from the cited classifications.
The symbol \(\underline c\) denotes their constant real line. All curvatures use
the Grassmannian normalization of
\eqref{eq:quadric-metric-normalization}.

\begin{table}[htbp]
\centering\small
\caption{Known low-dimensional families and their relation to the
finite-isotropy irreducible classification.}
\label{tab:known-low-dimensional}
\begin{tabular}{@{}llcp{0.39\textwidth}@{}}
\toprule
Target & Family & \(K\) & Relation to the hypotheses \\
\midrule
\(Q_2\) & \(\Phi_0^{(1)}\) & \(2\) & Reducible, strongly isotropic \\
 & \([UV_1^{(2)}]\oplus\underline c\) & \(1\) & Reducible, strongly isotropic \\
\midrule
\(Q_3\) & \(\Phi_0^{(4)}\) & \(1/2\) & Reducible, \(r=3\) \\
 & \(\Phi_1^{(4)}\) & \(1/5\) & Included: \(\ell=1,m=4,r=1\) \\
\midrule
\(Q_4\) & \(\Phi_0^{(3)}\) & \(2/3\) & Reducible, finite isotropy \\
 & \(\Phi_0^{(2)}\) & \(1\) & Reducible, finite isotropy \\
 & \(\Phi_0^{(2)}\) & \(1\) & Reducible, strongly isotropic \\
 & \([UV_2^{(4)}]\oplus\underline c\) & \(1/3\) & Reducible, strongly isotropic \\
 & \(\Phi_1^{(2)}\) & \(1/2\) & Irreducible, strongly isotropic \\
\midrule
\(Q_5\) & \(\Phi_1^{(6)}\) & \(1/8\) & Included: \(\ell=1,m=6,r=3\) \\
 & \(\Phi_2^{(6)}\) & \(1/11\) & Included: \(\ell=2,m=6,r=1\) \\
 & \(\Phi_0^{(m)},\ 3\leq m\leq6\) & \(2/m\) & Reducible, \(r=1\) \\
 & \(\Phi_0^{(6)}\) & \(1/3\) & Reducible, \(r=5\) \\
\bottomrule
\end{tabular}
\end{table}

The three finite-isotropy irreducible models in
Table~\ref{tab:known-low-dimensional} are the depth-one or depth-two
self-dual models of Example~\ref{ex:self-dual-family}; their parameters
give exactly the curvatures in \eqref{eq:main-curvature}. The remaining
families are reducible or strongly isotropic. For \(Q_2\), the BEW
bound \eqref{eq:intro-depth-bound} already excludes an irreducible
finite-isotropy tower, since it would require \(4\leq m\leq3\).

The nonhomogeneous \(Q_4\) example belongs to the reducible
\(\Phi_0^{(3)}\) family, and the noncongruent \(Q_5\) families in
\citep[Lemma~5]{jiao2019q5} belong to the reducible \(r=1\) row.
Their real Gram data distinguish quadric congruence classes even when
their projective representatives are Veronese. This distinction is
retained by Proposition~\ref{prop:real-gram}.

\Needspace{6\baselineskip}
\subsection{Examples}

\begin{example}[A self-dual family at every actual depth]
\label{ex:self-dual-family}
Fix \(\ell\geq1\) and an even integer \(m\geq2\ell+2\), and set
\(n=m+1\). With indices \(0\leq h,k\leq m\), put
\[
 W_{hk}=(-1)^h\delta_{h+k,m},\qquad
 P_\pm=\frac{I\pm W}{2},\qquad U=P_++\mathrm iP_-.
\]
The real symmetric involution \(W\) gives \(U^*U=I\) and
\(U^{\transp}U=W\). The binomial identity
\[
 (V_0^{(m)}(z))^{\transp}WV_0^{(m)}(t)=(t-z)^m
\]
shows that the finite-isotropy exponent is \(p=m\). The real Gram
matrix has rank \(m+1=n\). Corollary~\ref{cor:veronese-converse}
therefore applies with
\[
 s=\frac m2-\ell-1\geq0,\qquad r=m-2\ell-1,
 \qquad K=\frac{2}{(2\ell+1)m-2\ell^2}.
\]
\end{example}

Table~\ref{tab:veronese-examples} records the geometric parameters for
several members of this family.

\begin{table}[H]
\centering\small
\caption{Parameters in the self-dual Veronese family.}
\label{tab:veronese-examples}
\begin{tabular}{@{}cccccc@{}}
\toprule
\(\ell\) & \(m\) & \(n\) & \(r\) & \(s\) & \(K\) \\
\midrule
1 & 4 & 5 & 1 & 0 & \(1/5\) \\
2 & 6 & 7 & 1 & 0 & \(1/11\) \\
3 & 8 & 9 & 1 & 0 & \(1/19\) \\
1 & 6 & 7 & 3 & 1 & \(1/8\) \\
2 & 10 & 11 & 5 & 2 & \(1/21\) \\
3 & 14 & 15 & 7 & 3 & \(1/40\) \\
\bottomrule
\end{tabular}
\end{table}

\begin{example}[Real fullness with a smaller complex seed span]
Use the same square \(U\) and put
\[
 \widehat U=
 \begin{pmatrix}\frac35U\\[1mm]\frac{4\mathrm i}5U\end{pmatrix}
 :\mathbb C^{m+1}\hookrightarrow\mathbb C^{2m+2}.
\]
Then \(\widehat U^*\widehat U=I\) and
\(\widehat U^{\transp}\widehat U=-7W/25\). Its bilinear kernel
is \(-7(t-z)^m/25\). The real Gram matrix is positive definite,
with eigenvalues \(9/25\) and \(16/25\), each of multiplicity
\(m+1\). Thus the corresponding real layer is linearly full in
\(G(2,2m+2;\mathbb R)\), with the same depth, finite isotropy,
and curvature. In particular, real fullness does not imply
\(m=n-1\).
\end{example}

The theorem retains the geometric distinctions between immersion and
branching, finite and strong isotropy, and irreducible and reducible
sequences.

\appendix

    \section{Graph Connection and Contact Identities}
    \label{app:graph-connection}

    We give the full simultaneous-induction proof of
    Lemma~\ref{lem:simultaneous}, including its auxiliary derivative
    identities. All vector identities are read on the common monic
    regular set, with mixed index sums at most \(2\ell+1\).

    The simultaneous induction has seven stages:
    \begin{enumerate}
    \item compute the backward derivative of \(v_{j-1}\);
    \item define the forward direction \(u_{j-1}\) and prove it is
    orthogonal to the old graph;
    \item use the contact condition to restrict
    \(\bar\partial u_{j-1}\);
    \item impose induced holomorphicity to obtain the scalar
    replacement equation;
    \item establish orthogonality at the new level \(j\);
    \item derive the new column of the graph connection;
    \item propagate the contact condition to the new graph.
    \end{enumerate}
    This is the logical order in the proof: each stage uses only
    identities already established at the preceding level.

    \begin{proof}[Proof of Lemma~\ref{lem:simultaneous}]
    	For \(j=0\), \(v_0=w_0^0=f_0\) and \(S^0=(0)\), so
    	\eqref{eq:orthogonality-v-w} and \eqref{eq:S-connection} hold and
    	\eqref{eq:contact-condition} is vacuous. Assume the three assertions
    	at level \(j-1\), where \(j\geq1\).
    	
    	\smallskip
    	\noindent\emph{1. The backward derivative.}
    	Differentiating
    	\(v_{j-1}=\sum_{i=0}^{j-1}p_i^{j-1}w_i^{j-1}\) with the connection
    	formula at level \(j-1\), and using \(q_{j-1}^{j-1}=0\), gives
    	\begin{align}
    		 \bar\partial v_{j-1}
    		&=\sum_{i=0}^{j-2}d_i^{j-1}w_i^{j-1},
    		\label{eq:backward-v}\\
    		d_i^{j-1}
    		&= \bar\partial p_i^{j-1}
    		-\lambda_ip_{i+1}^{j-1}.
    	\end{align}
    	
    	\smallskip
    	\noindent\emph{2. The forward direction is orthogonal to the old graph.}
    	Set
    	\begin{equation}\label{eq:forward-direction}
    		u_{j-1}= \partial v_{j-1}
    		-\bigl(\partial\log\norm{v_{j-1}}^2\bigr)v_{j-1}.
    	\end{equation}
    	The induction hypothesis and \eqref{eq:backward-v} give
    	\(H(v_{j-1},\bar\partial v_{j-1})=0\), while the connection formula
    	at level \(j-1\) gives
    	\(H(v_{j-1},\bar\partial w_i^{j-1})=0\). Differentiating
    	\(H(v_{j-1},w_i^{j-1})
    	=\norm{v_{j-1}}^2\delta_{i,j-1}\) therefore yields
    	\begin{equation}\label{eq:U-orthogonal}
    		H(u_{j-1},w_i^{j-1})=0,
    		\qquad 0\leq i\leq j-1.
    	\end{equation}
    	
    	The Frenet support and \eqref{eq:graph-decomposition} show that
    	\begin{equation}\label{eq:U-isotropy}
    		\begin{aligned}
    		u_{j-1}&=f_j+r_{j-1},
    	\qquad r_{j-1}\in\overline{\mathbb L_{j-1}},\\
    		b(u_{j-1},u_{j-1})&
    		=b(u_{j-1},v_{j-1})
    		=b(u_{j-1},\overline{v_{j-1}})=0.
    		\end{aligned}
    	\end{equation}
    	
    	\smallskip
    	\noindent\emph{3. The contact condition restricts
    		\(u_{j-1}\).}
    	Differentiating the definition of \(u_{j-1}\) in
    	\eqref{eq:forward-direction}, we have
    	\begin{equation}\label{eq:U-backward-raw}
    		 \bar\partial u_{j-1}
    		=\partial(\bar\partial v_{j-1})
    		-(\bar\partial\partial\log\norm{v_{j-1}}^2)v_{j-1}
    		-(\partial\log\norm{v_{j-1}}^2)\bar\partial v_{j-1}.
    	\end{equation}
    	The Frenet equations and \eqref{eq:backward-v} show that this vector
    	lies in
    	\(\mathbb L_{j-1}\oplus_H\overline{\mathbb L_{j-1}}\).
    	
    	For \(i<j-1\), \(b(\bar\partial v_{j-1},w_i^{j-1})=0\) by
    	\(b\)-isotropy. Differentiating and using the induction hypothesis
    	gives
    	\begin{equation}
    		\begin{aligned}
    		b(\partial\bar\partial v_{j-1},w_i^{j-1})
    		&=-b(\bar\partial v_{j-1}, \partial w_i^{j-1})\\
    		&=-\sum_{k=0}^{j-2}d_k^{j-1}C_{ik}^{j-1}=0,\\
    		b(\bar\partial u_{j-1},w_i^{j-1})&=0.
    	\end{aligned}
    	\end{equation}
    	The remaining terms in \eqref{eq:U-backward-raw} lie in the
    	\(b\)-isotropic space \(\mathbb L_{j-1}\). Since
    	\begin{equation}
    		b(\overline{v_{j-1}},w_i^{j-1})
    		=\norm{v_{j-1}}^2\delta_{i,j-1}
    	\end{equation}
    	and the pairing is nondegenerate, there exist \(Z_{j-1}\in\mathbb L_{j-1}\)
    	and a scalar function \(c_{j-1}\) satisfying \eqref{eq:Ubar-decomposition} at level \(j-1\):
    	\begin{equation*}
    		 \bar\partial u_{j-1}=Z_{j-1}+c_{j-1}\overline{v_{j-1}}.
    	\end{equation*}
    	
    	\smallskip
    	\noindent\emph{4. The induced holomorphicity condition.}
    	By \eqref{eq:deviation-replacement},
    	\(v_j=u_{j-1}+x_{j-1}\overline{v_{j-1}}\), and \(v_j\) spans a
    	holomorphic line for the induced structure on
    	\(E=\operatorname{span}\{u_{j-1},\overline{v_{j-1}}\}\).
    	Using the conjugate of \eqref{eq:forward-direction},
    	\(\bar\partial\overline{v_{j-1}}
    	=\overline{u_{j-1}}
    	+(\bar\partial\log\norm{v_{j-1}}^2)\overline{v_{j-1}}\)
    	and \eqref{eq:Ubar-decomposition} gives
    	\begin{equation}\label{eq:v-next-derivative}
    		\begin{aligned}
    		 \bar\partial v_j
    		&=Z_{j-1}+x_{j-1}\overline{u_{j-1}}\\
    		&\quad+
    		\left(c_{j-1}+\bar\partial x_{j-1}
    		+x_{j-1}\bar\partial\log\norm{v_{j-1}}^2\right)
    		\overline{v_{j-1}}.
    		\end{aligned}
    	\end{equation}
    	By \eqref{eq:U-orthogonal} and \eqref{eq:U-isotropy},
    	\(Z_{j-1}+x_{j-1}\overline{u_{j-1}}\) is Hermitian orthogonal to
    	\(E\). Hence the induced holomorphicity condition forces
    	\eqref{eq:scalar-replacement}, namely
    	\begin{equation*}
    		c_{j-1}+\bar\partial x_{j-1}
    		+x_{j-1}\bar\partial\log\norm{v_{j-1}}^2=0.
    	\end{equation*}
    	Consequently,
    	\begin{equation}\label{eq:bar-v-next}
    		 \bar\partial v_j=Z_{j-1}+x_{j-1}\overline{u_{j-1}},
    		\qquad
    		b(\bar\partial v_j,w_i^{j-1})=0\quad(0\leq i<j).
    	\end{equation}
    	Here the pairing with \(\overline{u_{j-1}}\) vanishes by
    	\eqref{eq:U-orthogonal}, and \(Z_{j-1}\in\mathbb L_{j-1}\).
    	
    	\smallskip
    	\noindent\emph{5. The orthogonality relation at the new level.}
    	Since \(b(v_j,v_j)=0\) by \eqref{eq:U-isotropy}, the
    	graph-extension construction of Definition~\ref{def:S-and-w} and
    	Lemma~\ref{lem:graph-properties} gives \(v_j\in\mathbb L_j\). Moreover,
    	\(u_{j-1}\) and \(\overline{v_{j-1}}\) are Hermitian orthogonal to
    	\(\mathbb L_{j-1}\), and
    	\(H(v_j,\overline{f_j})=b(v_j,f_j)=0\). The update formula
    	\eqref{eq:w-update} and the graph representation therefore yield
    	\begin{equation}
    		H(v_j,w_i^j)=\norm{v_j}^2\delta_{ij},
    		\qquad 0\leq i\leq j.
    	\end{equation}
    	This proves assertion (1) at level \(j\).
    	
    	\smallskip
    	\noindent\emph{6. The new column of the connection equation.}
    	For \(0\leq i<j\), differentiating
    	\(t_i^j=b(v_j,w_i^{j-1})\), substituting
    	\eqref{eq:bar-v-next} and the level-\((j-1)\) connection formula,
    	and using \(p_j^j=1\) and
    	\(b(v_j,\overline{f_j})=\norm{f_j}^2\), we obtain
    	\begin{align}
    		 \bar\partial t_i^j
    		&=-\lambda_{i-1}
    		b(v_j,w_{i-1}^{j-1})
    		+\frac{S_{j-1,i}^{j-1}}{\norm{f_{j-1}}^2}
    		b(v_j,\overline{f_j})\\
    		&=-\lambda_{i-1}t_{i-1}^j
    		-\lambda_{j-1}S_{i,j-1}^{j-1}.
    	\end{align}
    	Since \(S_{ij}^j=t_i^j\) and
    	\(S_{j-1,i}^{j-1}=-S_{i,j-1}^{j-1}\), this gives
    	\begin{equation}\label{eq:new-column-connection}
    		\boxed{
    			 \bar\partial S_{ij}^{j}
    			=-\lambda_{i-1}S_{i-1,j}^{j}
    			-\lambda_{j-1}S_{i,j-1}^{j-1},
    			\qquad 0\leq i<j.}
    	\end{equation}
    	The old principal block is unchanged, so all connection equations
    	involving only old indices remain valid. Together with skew-symmetry,
    	\eqref{eq:new-column-connection} supplies all remaining entries of
    	\eqref{eq:S-connection}. This proves assertion (2).
    	
    	\smallskip
    	\noindent\emph{7. Propagation of the contact condition.}
    	Since \(b(w_i^j,w_k^j)=0\), applying \(\partial\) gives
    	\begin{equation}\label{eq:C-skew}
    		C_{ik}^j=-C_{ki}^j.
    	\end{equation}
    	For \(i<j-1\), we have
    	\begin{equation}
    		\begin{aligned}
    		t_i^j
    		&=b(v_j,w_i^{j-1})
    		=b(u_{j-1},w_i^{j-1})\\
    		&=b(\partial v_{j-1},w_i^{j-1})\\
    		&=-b(v_{j-1}, \partial w_i^{j-1})\\
    		&=-\sum_{k=0}^{j-1}p_k^{j-1}C_{ik}^{j-1}\\
    		&=-C_{i,j-1}^{j-1}
    		=C_{j-1,i}^{j-1}.
    		\end{aligned}
    	\end{equation}
    	Here we used \eqref{eq:forward-direction}, the induction hypothesis,
    	\(b(\overline{v_{j-1}},w_i^{j-1})=0\), and
    	\(p_{j-1}^{j-1}=1\). This proves
    	\eqref{eq:contact-new-column}; its
    	bilinear-pairing form follows from the definition of \(C\).
    	
    	By \eqref{eq:w-update},
    	\begin{gather}
    		w_i^j=w_i^{j-1}
    		-\frac{t_i^j}{\norm{f_j}^2}\overline{f_j},
    		\qquad i<j,\\
    		\partial\overline{f_j}
    		=-\lambda_{j-1}\overline{f_{j-1}},
    	\end{gather}
    	because \(b(\overline{f_j},w_k^{j-1})=0\). Substitution into
    	\(C_{ik}^j=b(\partial w_i^j,w_k^j)\) gives, for \(i,k<j\),
    	\begin{equation}\label{eq:contact-update}
    		C_{ik}^j
    		=C_{ik}^{j-1}
    		-\delta_{i,j-1}t_k^j+\delta_{k,j-1}t_i^j.
    	\end{equation}
    	If \(i,k<j-1\), the right-hand side vanishes by induction. If
    	\(i<j-1\) and \(k=j-1\), then
    	\begin{equation}
    		C_{i,j-1}^j
    	=C_{i,j-1}^{j-1}+t_i^j
    	=-t_i^j+t_i^j=0,
    	\end{equation}
    		where we used \eqref{eq:contact-new-column} and
    	\eqref{eq:C-skew}. The remaining cases follow from
    	\eqref{eq:C-skew}, so
    	\begin{equation}
    		C_{ik}^j=0,\qquad 0\leq i,k<j.
    	\end{equation}
    	This proves assertion (3) and completes the induction.

    	The terminal level requires no extra replacement. Repeating
    	steps 1--3 with the established level-\(\ell\) graph gives
    	\eqref{eq:terminal-graph-connection}, namely
    	\begin{equation*}
    	 \bar\partial v_\ell\in
    	 \operatorname{span}\{w_0^\ell,\ldots,w_{\ell-1}^\ell\},
    	 \qquad
    	 \bar\partial u_\ell=Z_\ell+c_\ell\overline{v_\ell},
    	 \quad Z_\ell\in\mathbb L_\ell,
    	\end{equation*}
    	where \(u_\ell=\partial v_\ell-
    	(\partial\log\norm{v_\ell}^2)v_\ell\). The same computation
    	also gives \(u_\ell\perp_H\mathbb L_\ell\) and
    	\(b(u_\ell,v_\ell)=0\), completing the terminal-level identities.
    \end{proof}

\section{The Bordered-Pfaffian Identity}
\label{app:bordered-pfaffian}

We give the block expansion and sign calculation underlying
\eqref{eq:rigidity-pf-border}. This also fixes the sign of the
mixed kernel in Lemma~\ref{lem:rigidity-mixed}.

Let \(K\) be a complex skew-symmetric matrix of order \(d\), and
put \(M=I+KK^*\). Use \(\operatorname{Pf}(K_\varnothing)=1\)
and set
\[
 \mathscr N(K)=\sum_{I\ \mathrm{even}}|\operatorname{Pf}(K_I)|^2,
 \qquad
 \mathcal B=\begin{pmatrix}K&I\\-I&-\overline K\end{pmatrix}.
\]
All principal submatrices use increasing index order. We claim that
\begin{equation}\label{eq:pfaffian-block-expansion}
 \operatorname{Pf}\mathcal B=(-1)^{d(d-1)/2}\mathscr N(K).
\end{equation}
Indeed, label the two blocks by \(e_1,\ldots,e_d\) and
\(g_1,\ldots,g_d\). In any nonzero matching, the indices
paired within the first block form an even set \(I\); the
identity off-diagonal block forces the same set to be paired
within the second block. Write \(|I|=2h\) and
\(J=\{1,\ldots,d\}\setminus I\). Reordering to
\((e_I,e_J,g_J,g_I)\) gives two shuffle signs whose product is
\((-1)^{|I||J|}=1\). The internal lower-block matching contributes
\[
 \operatorname{Pf}(-\overline K_I)
 =(-1)^h\overline{\operatorname{Pf}(K_I)},
\]
whereas the remaining identity matching has sign
\((-1)^{(d-2h)(d-2h-1)/2}\). Their combined sign is
\[
 (-1)^{h+(d-2h)(d-2h-1)/2}=(-1)^{d(d-1)/2}.
\]
Summing proves \eqref{eq:pfaffian-block-expansion}. In particular,
\(\mathscr N(K)\geq1\) makes \(\mathcal B\) invertible.

Append two indices after the lower block and form
\[
 \mathcal A=
 \begin{pmatrix}
 K&I&\mathsf v&\mathsf t\\
 -I&-\overline K&0&0\\
 -\mathsf v^{\transp}&0&0&\omega\\
 -\mathsf t^{\transp}&0&-\omega&0
 \end{pmatrix}.
\]
In its matching expansion, the two new indices can pair only
with each other or with upper indices. For a fixed lower internal
set \(I\), the upper internal part is therefore the bordered
matrix in \eqref{eq:rigidity-pf-border}. Moving the two new
indices next to \(e_I\) has positive sign, since they pass
every intervening index as a block of size two. The same shuffle
calculation consequently gives
\begin{equation}\label{eq:pfaffian-bordered-expansion}
 \operatorname{Pf}\mathcal A
 =(-1)^{d(d-1)/2}
 \sum_{I\ \mathrm{even}}
 \operatorname{Pf}
 \begin{pmatrix}
 K_I&\mathsf v_I&\mathsf t_I\\
 -\mathsf v_I^{\transp}&0&\omega\\
 -\mathsf t_I^{\transp}&-\omega&0
 \end{pmatrix}
 \overline{\operatorname{Pf}(K_I)}.
\end{equation}

It remains to evaluate the same Pfaffian by a congruence. Solving
\(\mathcal B(x,y)^{\transp}=(u,0)^{\transp}\) gives
\(x=-\overline K y\) and
\((I-K\overline K)y=u\). Since \(K^*=-\overline K\),
the upper-left inverse block is
\begin{equation}\label{eq:pfaffian-inverse-block}
 (\mathcal B^{-1})_{11}=-\overline K M^{-1}.
\end{equation}
Write
\[
 \mathcal E=\begin{pmatrix}\mathsf v&\mathsf t\\0&0\end{pmatrix},
 \qquad
 \mathcal D=\begin{pmatrix}0&\omega\\-\omega&0\end{pmatrix}.
\]
For an invertible even-order skew matrix \(A\), an even-order
skew matrix \(D\), and a rectangular matrix \(E\), put
\[
 R=\begin{pmatrix}I&-A^{-1}E\\0&I\end{pmatrix}.
\]
Since \(A^{-1}\) is skew-symmetric, direct block multiplication gives
\[
 R^{\transp}
 \begin{pmatrix}A&E\\-E^{\transp}&D\end{pmatrix}R
 =
 \begin{pmatrix}A&0\\0&D+E^{\transp}A^{-1}E\end{pmatrix},
 \qquad\det R=1.
\]
The Pfaffian congruence and block-diagonal identities therefore yield
\begin{equation}\label{eq:app-pfaffian-schur}
 \operatorname{Pf}
 \begin{pmatrix}A&E\\-E^{\transp}&D\end{pmatrix}
 =
 \operatorname{Pf}(A)\,
 \operatorname{Pf}(D+E^{\transp}A^{-1}E).
\end{equation}
Apply \eqref{eq:app-pfaffian-schur} with
\(A=\mathcal B\), \(D=\mathcal D\), and \(E=\mathcal E\). By
\eqref{eq:pfaffian-inverse-block} and the skew-symmetry of
\(\mathcal B^{-1}\), the lower-right block is
\[
 \begin{pmatrix}
 0&\omega-\mathsf v^{\transp}\overline K M^{-1}\mathsf t\\
 -\omega+\mathsf v^{\transp}\overline K M^{-1}\mathsf t&0
 \end{pmatrix},
\]
whose Pfaffian is
\(\omega-\mathsf v^{\transp}\overline K M^{-1}\mathsf t\).
Consequently,
\[
 \operatorname{Pf}\mathcal A
 =\operatorname{Pf}\mathcal B\,
   (\omega-\mathsf v^{\transp}\overline K M^{-1}\mathsf t).
\]
Substituting \eqref{eq:pfaffian-block-expansion} and
\eqref{eq:pfaffian-bordered-expansion} cancels the common sign
and proves \eqref{eq:rigidity-pf-border}, including the empty
case \(d=0\), where \(\mathscr N(K)=1\) and
\(\operatorname{Pf}\mathcal A=\omega\).

\newpage
\bibliographystyle{plain}
\bibliography{02}

\Needspace{9\baselineskip}
\section*{AI Assistance Statement}

The final logical structure of the proof was determined entirely by the
authors. AI's role in this process was limited to iterative suggestions
and draft generation. The authors verified and corrected every step
individually and can independently justify each one. In addition, AI
assisted with the writing and revision of this paper.

\end{document}